\magnification =\magstep1
\baselineskip =13pt
\overfullrule =0pt
\ifnum\pageno=1\nopagenumbers

\input amstex
\documentstyle{amsppt}
\nologo

\vsize=500pt

\topmatter

\title{Derived Hilbert schemes }\endtitle
\rightheadtext{Derived Hilbert schemes}
\author{Ionu\c t Ciocan-Fontanine and
Mikhail M. Kapranov}\endauthor

\address{Department of Mathematics,
Northwestern University, Evanston, IL 60208}\endaddress

\email ciocan\@math.nwu.edu,  kapranov\@math.nwu.edu \endemail

\endtopmatter

\document

\def\op{\operatorname}
\def\lra{\longrightarrow}
\def\ind{{\lim\limits_{\longrightarrow}}}

\heading {Introduction}\endheading

\noindent {\bf (0.1)}  
The Derived Deformation Theory (DDT) program
(see \cite{Kon},  \cite{CK} for more details and historical references) seeks
to avoid the difficulties related to the singular nature of
the moduli spaces in geometry by ``passing to the derived category",
i.e., developing an appropriate version of the (nonabelian) derived
functor of the functor of forming the moduli space. The resulting
geometric objects are sought to  be not ordinary varieties or schemes but
rather
dg-schemes, i.e., geometric objects whose algebras of functions
are commutative differential graded (dg-)algebras and which  are considered
up to quasi-isomorphism. Moreover, they are expected to
be smooth in the sense that the corresponding dg-algebras
can be obtained from a smooth commutative algebra in the usual
sense by adding projective modules of generators in each
negative degree.

 At the present time, one can say that the DDT
program is well established in the formal case (i.e., one has
a good theory of the derived version of the formal neighborhood
of any point in the usual deformation space,  see \cite {BK} \cite{Hi1}
\cite{Man}). However, global derived moduli spaces
(of great interest for enumerative applications)
are much less  understood. 
 In \cite{CK}, we have constructed the derived
version of the first global algebro-geometric moduli space,
namely of the Grothendieck's $Quot$ scheme. The aim of the
present paper is to give a construction
(using a similar but different approach) of the derived Hilbert scheme. 
  
\vskip .3cm

\noindent {\bf (0.2)} While in the usual algebraic geometry,
the Hilbert scheme is a particular case of the $Quot$ schemes,
the two constructions  diverge in the derived
world. To see the difference, consider a smooth projective
variety $X$ and a  subvariety $Z\i X$
with Hilbert polynomial $h$ and sheaf of ideals
$\Cal I_Z$. Then $Z$
represents a point $[Z]$ in $Hilb_h(X) = Quot_h(\Cal O_X)$.
The tangent space to the dg-manifold $RQuot_h(\Cal O_X)$  at the point $[Z]$
is a cochain complex with  cohomology
$$H^iT^\bullet_{[Z]} RQuot_h(\Cal O_X) = \op{Ext}^{i}_{\Cal O_X}
(\Cal I_Z, \Cal O_Z).\leqno (0.2.1)$$
In the present paper we construct a dg-manifold $RHilb_h^{LCI}(X)$
whose degree 0  truncation $\pi_0 RHilb_h^{LCI}(X)$ is the
open part in $Hilb_h(X)$ consisting of locally complete intersections
and whose tangent dg-space at a point $[Z]$ for $Z$ as above
is given by
$$H^iT^\bullet_{[Z]}RHilb_h^{LCI}(X) = H^i(Z, N_{Z/X}),\leqno (0.2.2)$$
where $N_{Z/X}$ is the normal bundle. 
See (4.2) for precise statements
in the case of arbitrary $Z$, when the relative cotangent
complex is involved. 
When $h$ is identically 1, we get $RHilb_h^{LCI}(x) = Hilb_h(X) = X$,
but $RQuot_h(\Cal O_X)$ is not identified with $X$ in any way,
because for $Z=\{x\}$ the tangent dg-space at $x$ to $RQuot$
has by (0.2.1) the $i$th cohomology equal to $\Lambda^{i+1}(T_xX)$.

The  dg-manifolds $RQuot$ from \cite{CK} 
 are suitable
for the construction of the derived moduli spaces
of vector bundles on a fixed algebraic variety, as explained in
\cite{CK}, Remark 4.3.8. In contrast, the dg-schemes $RHilb$
constructed here, should play the same role for the
derived moduli spaces of algebraic varieties themselves.
For example, we use them to construct two types
of geometric derived moduli spaces:

\vskip .2cm

\noindent $\bullet$ The derived space of maps 
  $\op{RMap}(C, Y)$ from a fixed projective scheme $C$
to a fixed smooth projective variety $Y$. See (4.3.6)

\vskip .2cm

\noindent $\bullet$ The derived stack of stable
degree $d$ maps $R\overline {M}_{g,n}
(Y, d)$ from $n$-pointed nodal curves to a given smooth projective
variety $Y$,  see Section 5. The non-derived stacks
of stable maps $\overline{M}_{g,n}(Y, d)$ were introduced
by Kontsevich \cite{Kon}, see also \cite{BM} \cite{FP}, and
turned out to be extremely important in the
mathematical theory of Gromov-Witten invariants.
However, Kontsevich in fact proposed in {\it loc cit.} to construct
the derived version of $\overline{M}_{g,n}(Y, d)$ as well.
Our work  carries out this proposal. 

\vskip .3cm

\noindent {\bf (0.3)} Similarly to the case of
 the usual algebro-geometric
moduli spaces, it would be nice to characterize
$RHilb$ and $RQuot$ in terms of representability of
some functors. This is not easy, however, as the functors
should be considered on the {\it derived} category of
dg-schemes (with quasi-isomorphisms inverted)
and for morphisms in this localized category there is 
currently no explicit description. This issue should
be probably addressed in a wider foundational context
for dg-algebraic geometry, in which one should
systematically work with objects glued from dg-schemes
in our present sense by means of gluing maps which
are only quasi-isomorphisms on  pairwise intersections,
satisfying cocycle conditions only up to homotopy on
triple intersections etc.

Nevertheless, 
our derived Hilbert scheme $RHilb$
carries a certain natural family of commutative
dg-algebras which over the usual Hilbert scheme
is just the family of truncations of the 
graded coordinate rings of $Z\i X$. 
In contrast, the dg-scheme
 $RQuot$ carries a family of dg-modules, not algebras. 

Accordingly,
the main technical step in constructing $RHilb$
is the construction of the derived moduli space $RCA(W)$ of structures
of a commutative algebra on a given finite-dimensional vector
space $W$, see \S 3. 
We achieve this by using a certain dg-resolution of the operad $\Cal Com$
describing commutative algebras. Derived moduli spaces of operad
algebras were the subject of several recent papers \cite {Re} \cite{Hi2}
\cite{KS} although our approach differs from the approaches
of these papers. 
The similar step 
in the construction of $RQuot$ in \cite{CK} involved
the derived space of structures of an $A$-module on $W$, where
$A, W$ are fixed algebra and vector space, respectively. That
was achieved by using a dg-resolution of $A$.

\vskip .3cm

\noindent {\bf (0.4)} Let us now describe the
content of the paper. In Section 1 we give an identification of the
Hilbert scheme of a projective scheme $X$ with graded coordinate
ring $A$ with the scheme of graded ideals of the finite-dimensional
graded algebra (without unit) $A_{[p,q]}$ for $0\ll p\ll q$.
This result seems to be new and interesting by itself. It improves
Theorem 1.4.1 of \cite{CK}. 

Section 2 provides some background information on the cotangent complexes
and Harrison homology.
The only possibly new result here is Theorem 2.2.2 which connects
the relative cotangent complex of a morphism of graded algebras
with the relative cotangent complex of the morphism of their
projective spectra.

In Section 3 we develop the central construction of the paper:
the derived scheme of ideals in a finite-dimensional commutative
algebra. This is parallel to (but different from) \cite{CK}
where the corresponding role was played by the
derived schemes of $A$-submodules in a finite-dimensional
module over a given algebra $A$. 

In Section 4 we apply the algebraic formalism of Section 3
to the truncated graded coordinate rings $A_{[p,q]}$, give
the constructions of the dg-manifolds $RHilb_h^{\leq m}(X)$ and
$RHilb_h^{LCI}(X)$ and prove the properties stated in (0.3).
We also construct the derived space of maps $\op{RMap}(C, Y)$.

Finally, in Section 5 we relativize the above construction
to the case of a base stack $\Cal S$ of finite type.
In particular, we develop  the beginnings of a theory of dg-stacks. 
We apply this to the construction of the derived
stack of stable maps. 

\vskip .3cm

\noindent {\bf (0.5)} We would like to thank V. Drinfeld and K. Behrend 
who made us realize that there should be not one, but two derived versions
of the Hilbert scheme. 
Both authors were partially supported by NSF grants.

\vfill\eject

\heading {1.Hilbert scheme as the variety of ideals}\endheading

\noindent {\bf (1.1) Generalities on $Hilb$ and $Quot$.}
We work over a fixed field $\Bbb K$. 
Let $X\i \Bbb P^n$
 be a  projective scheme over $\Bbb K$ and $\Cal O_X(1)$
be the very ample sheaf defining the projective embedding. For a coherent
sheaf $\Cal G$ on $X$ we denote $\Cal G(i)=\Cal G\otimes\Cal O_X(i)$ and 
denote by $h^{\Cal G}(t)\in \Bbb Q[t]$ the Hilbert polynomial of $\Cal G$.
Thus $h^{\Cal G}(i) = \chi (\Cal G(i))$ for $i\in \Bbb Z$. 

Fix a polynomial $h\in \Bbb Q[t]$ and let $Hilb_h(X)$
be the Hilbert scheme parametrizing closed subschemes $Z\subset X$ such that
$h^{\Cal O_Z}=h$. For such $Z$ we denote by $\Cal I_Z\i \Cal O_X$
the sheaf of ideals of $Z$ and by $[Z]$ the $\Bbb K$-point of
$Hilb_h(X)$ corresponding to $Z$. Note that the full Hilbert scheme
of $X$ does not depend on the projective embedding, but
the way it is decomposed into the disjoint union of the 
subschemes $Hilb_h(X)$,
does. Should we need to emphasize the dependence of these subschemes
on the embedding, we will use the notation $Hilb_h(X, \Bbb P^n)$. 

More generally, if $\Cal F$ is a coherent sheaf on $X$, we have the $Quot$
scheme $Quot_h(\Cal F)$ parametrizing flat families of quotients
$\Cal F\twoheadrightarrow \Cal G$ with
$h^{\Cal G}=h$. The Hilbert scheme is a particular case corresponding
to $\Cal F = \Cal O_X$. See \cite{Gr} \cite{Kol} \cite{Vi} for more background. 

Recall Grothendieck's Grassmannian embedding of the Quot scheme, see {\it loc. cit.}
Set $k=h^{\Cal F}-h\in \Bbb Q[t]$. We will identify the polynomial
$k$ with the sequence of its values $k_i = k(i)$, $i\in \Bbb Z$.  
Fix $r\gg 0$. Associating to a quotient
$\Cal F\twoheadrightarrow \Cal G$ the linear subspace
$\text{Ker}\{H^0(X, \Cal F(r))\to H^0(X, \Cal G(r))\}$ defines
a regular map 
$$\alpha_r:
Quot_h(\Cal F)\to G(k(r), H^0(X, \Cal F(r))).$$ 

\proclaim {(1.1.1) Theorem} For $r\gg 0$ the map $\alpha_r$
identifies $Quot_h(\Cal F)$ with a closed subscheme of the
Grassmannian  $G(k(r), H^0(X, \Cal F(r)))$.

\endproclaim

\vskip .3cm

\noindent {\bf (1.2) The $A$-Grassmannian and the scheme of ideals.}
Let $A$ be an associative algebra over $\Bbb K$ (possibly
without unit) and $M$ a  left $A$-module
which is finite-dimensional (over $\Bbb K$).
Then we have the $A$-Grassmannian $G_A(k,M)\i G(k,M)$,
parametrizing $k$-dimensional submodules $V\i M$, see  n.(1.3) of \cite{CK}. 
More generally, let $A$ be a graded algebra and $M$ be a finite-dimensional
graded $A$-module. Then, for every sequence $k=(k_i)_{i\in {\Bbb Z}}$
of non-negative integers we have the graded $A$-Grassmannian
$G_A(k, M)$ which parametrizes graded submodules $V\i M$ with
$\dim_{\Bbb K}(V_i) = k_i$. 

Suppose now that $A$ is finite-dimensional. Then
we can take $M=A$; the scheme $G_A(k, A)$ will be denoted
$J(k, A)$ and called the {\it scheme of (left) ideals} in $A$
of dimension $k$. For such an ideal $I\i A$ we will denote
$[I]\in J(k, A)$ the corresponding $\Bbb K$-point. 

Suppose, moreover, that $A$ is commutative. Then left
ideals are two-sided and for such an ideal $I$ the quotient
$A/I$ is again a commutative algebra. In this case the Zariski tangent
spaces to $J(k, A)$ are given by:
$$T_{[I]}J(k, A) = \text{Hom}_A(I, A/I) = \text{Hom}_{A/I}(I/I^2, A/I).
\leqno (1.2.1)$$
Here the first equality is the general property of $A$-Grassmannians,
see n.(1.3.2) of  \cite{CK} and the second equality is due to the commutativity of $A$. 

Let now $A$ be a graded finite-dimensional commutative algebra and
$k=(k_i)$ be a sequence of integers, as before. We have then
the scheme $J(k, A)$ (the scheme of graded ideals)
parametrizing graded ideals $I\i A$ with $\dim(I_i) = k_i$. 
The description of the tangent space
in this case is modified as follows:
$$T_{[I]}J(k, A) =  \text{Hom}^0_{A/I}(I/I^2, A/I),
\leqno (1.2.2)$$
where $\text{Hom}^0$ means the space of homomorphisms of degree 0. 

\vskip .3cm

\noindent {\bf (1.3) $Hilb$ as the scheme of ideals.}
We return to the situation of (1.1) and let $A=\bigoplus_i H^0(X, \Cal O_X(i))$
be the homogeneous coordinate ring of $X$. For $p\geq 0$ let
$A_{\geq p} = \bigoplus_{i\geq p} H^0(X, \Cal O_X(i))$ be the
truncation of $A$ in degrees $\geq p$. This is an ideal in $A$;
we denote by $A_{[p,q]} = A_{\geq q}/A_{\geq p}$, $p\leq q$,
the finite-dimensional truncation in degrees from $p$ to $q$.
Both $A_{\geq p}$ and $A_{[p,q]}$ are graded commutative algebras
without unit. In particular, the construction of the graded
scheme of ideals is applicable to $A_{[p,q]}$. 

Suppose that $p\gg 0$. Then  for any $q\geq p$
the Grassmannian embeddings $\alpha_r$, $p\leq r\leq q$, from (1.1) define a morphism
$\beta_{[p,q]}: Hilb_h(X)\to J(k, A_{[p,q]})$. 
Here $k$ is the sequence $(k_p = h^{\Cal O_X}(p)-h(p), ..., 
k_q = h^{\Cal O_X}(q)-h(q))$. 

Let us now formulate the main result of this section.

\proclaim{(1.3.1) Theorem} If $0\ll p\ll q$, then the
morphism $\beta_{[p,q]}$ identifies $Hilb_h(X)$
with the scheme of graded ideals $J(k, A_{[p,q]})$.
\endproclaim

To begin the proof, recall a 
particular case of a result proved in \cite{CK} (1.4.1).

\proclaim {(1.3.2) Theorem} (a) For $0\ll p\ll q$ the natural
embedding $\alpha_{[p,q]}: Hilb_h(X)\to G_A(k, A_{[p,q]})$ is
an isomorphism. 

(b)  More precisely, (a) is true whenever $p$ is
such that the map $\alpha_p$ from (1.1.1) is an embedding and
$q$ is large enough compared with $p$. 
\endproclaim

Note that we have an embedding
$$\phi: G_A(k, A_{[p,q]}) \hookrightarrow J(k, A_{[p,q]})
\leqno (1.3.3)$$
since an $A$-submodule is automatically an $A_{[p,q]}$-submodule.
So Theorem 1.3.1 strengthens Theorem 1.3.2 which we will nevertheless
use as a starting point. Also, in the proof of Theorem 1.3.1
we can and will assume that $\Bbb K$ is algebraically closed. 

Choose $p>0$  such that $\alpha_p$ is an embedding and consider the
Veronese subalgebra $B=\bigoplus_j A_{pj}$, graded so that
$B_j = A_{pj}$. Let $\Bbb P^n = \Bbb P(A_1^*)$ be the ambient space of the 
initial projective
embedding of $X$. Then $\Bbb P^N = \Bbb P(A_p^*) = \Bbb P(B_1^*)$
is the ambient space of another projective embedding
known as the $p$-fold Veronese embedding. In other words,
$X=\text{Proj}(A) = \text{Proj}(B)$.
Note that
on $X$ we have $\Cal O_{\Bbb P^N}(1) \simeq \Cal O_{\Bbb P^n}(p)$. 
The independence of the full Hilbert scheme on the
projective embedding implies that
$$Hilb_h(X, \Bbb P^n) = Hilb_{h'}(X, \Bbb P^N), \quad 
\text {where}\quad h'(t) = h(pt).\leqno (1.3.4)$$
Let $r\geq 1$. Then, 
setting  $k'_i = k_{pi}$ for $i=1, ..., r$, we find
$$J(k', B_{[1,r]}) = G_B(k', B_{[1,r]}).\leqno (1.3.5)$$
By Theorem 1.3.2(b) we have that for $r\gg 0$
each of the schemes in (1.3.5) is identified with $Hilb_{h'}(X, \Bbb P^N)$,
i.e., with $Hilb_{h}(X, \Bbb P^n)$. 

We now take $q=pr, r\gg 0$ and construct an inverse to the
embedding (1.3.3). Given a graded ideal $I\i A_{[p, q]}$, we
form the ideal $J = \bigoplus_j I_{pj}$ in $B_{[1,r]}$
obtained by taking every $p$th graded component of $I$. 
Let $\Cal J\i \Cal O_X$ be the sheaf of ideals associated to
$J$ by (1.3.4) and (1.3.5). By construction, $\Cal O_X/\Cal J$
has  Hilbert polynomial $h'$ with respect to the embedding
of $X$ into $\Bbb P^N$, and therefore  has  Hilbert
polynomial $h$ with respect to the embedding of $X$ into
$P^n$, so it represents a point of $Hilb_h(X, \Bbb P^n)$.
Define now
$$\psi(I) = \bigoplus_{i=p}^q H^0(X, \Cal J\otimes\Cal O_{\Bbb P^n}(i)).$$
 
Our theorem follows from the next lemma.

\proclaim {(1.3.6) Lemma} (a) $\psi(I)\in G_A(k, A_{[p,q]})$.

(b) The correspondence $I\mapsto \psi(I)$ gives rise to
a morphism of schemes $\psi:  J(k, A_{[p,q]})
\to  G_A(k, A_{[p,q]})$.

(c) The morphisms $\phi$ and $\psi$ are mutually inverse. 

\endproclaim

\noindent {\sl Proof:}
(a) It is clear that $\psi(I)$ is a graded ideal in $A_{[p,q]}$. 
Because of the choice of $p$, the Hilbert function of the ideal
$\psi(I)$ is given by its Hilbert polynomial. (Among other things,
the choice of $p$ implies this property for all ideals from
$Hilb_h(X)$). In other words, $\dim_{\Bbb K} \psi(I)_j = k_j$, as claimed.

(b) First of all, passing from $I$ to $J$ is clearly a morphism
of schemes $J(k, A_{[p,q]})\to J(k', B_{[1,r]})$. The target of this
morphism is identified with $Hilb_h(X, \Bbb P^n)$. It remains
to notice that (again, because of the choice of $p$)
for any $i\geq p$ the correspondence
$$\Cal J \mapsto H^0(X, \Cal J\otimes \Cal O_{\Bbb P^n}(i))$$
is  a regular map (the Grothendieck embedding) of $Hilb_h(X, \Bbb P^n)$
into a Grassmannian.

\vskip .1cm

 For (c), it suffices to prove that  $\phi\psi=\text{Id}$. We
do this at the level of $\Bbb K$-points; the proof for $S$-points,
where $S$ is any scheme is similar.
So we start with a graded ideal
$I\i A_{[p,q]}$ with $\dim(I_i)=k_i$. For any $\nu =0,1, ..., p-1$ 
consider the
graded $B$-module
$$I^\nu = \bigoplus_{i\equiv \nu \,\,\text{mod}\, p} I_i.$$
By Theorem 1.3.2 applied to the $B$-module $\bigoplus_{i\equiv \nu \,\,
\text{mod}\, p} A_i$, we
associate to $I^\nu$ a coherent subsheaf $\Cal J^\nu 
\i \Cal O_{\Bbb P^n}(\nu)|_X$.
Our statement would follow from the identification 
$$\Cal J^\nu =
\Cal J\otimes\Cal O_{\Bbb P^n}(\nu)\leqno (1.3.7)$$
which we now prove. 
Fix some $\nu =0,1, ..., p-1$ and consider another
Veronese subalgebra
$$C = \bigoplus_j A_{j(p+\nu)}\i A, \quad C_j = A_{j(p+\nu)},
\leqno (1.3.8)$$
and the corresponding projective embedding $X\i \Bbb P^M = \Bbb P(C_1^*)$. 
Associate to $I$ the graded $C$-module
$$K = \bigoplus_j I_{j(p+\nu)},$$
where the summation is over such $j$ that $j(p+\nu)\leq q$. 
As before, $K$ gives rise to a sheaf of ideals $\Cal K\i \Cal O_X$
whose Hilbert polynomial with respect to the embedding $X\i \Bbb P^M$
is equal to $k''(t) = k((p+\nu)t)$ and thus its Hilbert polynomial with 
respect to $X\i \Bbb P^n$ is equal to $k$. 
Thus we have three $\Bbb K$-points of $Hilb_h(X, \Bbb P^n)$
represented by   $\Cal J$, $\Cal J^\nu \otimes\Cal O_{\Bbb P^n}(-\nu)$
 and $\Cal K$. We claim that they all coincide. 
For this, we consider their images under the Grothendieck
embedding into the Grassmannian of subspaces in 
$H^0(X, \Cal O_{\Bbb P^n}(p(p+\nu)))$ of dimension $k_{p(p+\nu)}$. 
Each of these images is equal to the subspace $I_{p(p+\nu)}$.
Since the Grothendieck embedding is indeed an embedding,
we find that the sheaves of ideals coincide, in particular, that $\Cal J = 
\Cal J^\nu \otimes\Cal O_{\Bbb P^n}(-\nu)$. This proves (1.3.7)
and hence Theorem 1.3.1 is proved. 

\vskip .3cm

\noindent {\bf (1.4) Relative version.} Let $S$ be a scheme of finite type
and $\Cal A$ be a commutative $\Cal O_S$-Algebra which, as an
$\Cal O_S$-module, is locally free of finite rank. We have then, for
each $k\geq 0$, the relative scheme of ideals $J(k, \Cal A/S)$
which is a closed $S$-subscheme in the relative Grassmannian
$G(k, \Cal A/S)$ such that for any $S$-scheme $a: Y\to S$ the
set $\op{Hom}_S(Y, J(k, \Cal A/S))$ is identified
with the set of sheaves of ideals $\Cal I\i q^*\Cal A$ on $Y$
which, as $\Cal O_Y$-modules,
 are locally free  direct summands of rank $k$. Similarly
for the case of a graded $\Cal O_S$-Algebra $\Cal A = \bigoplus_{i=a}^b
\Cal A_i$ and a sequence $k=(k_i)_{i=a}^b$ of
nonnegative integers. 

\vskip .1cm

Next, if $p: X\to S$ is a flat projective morphism and $\Cal O(1)$
is a fixed relative very ample sheaf, then for each $h\in \Bbb Q[t]$
we have the relative Hilbert scheme $Hilb_h(X/S)\to S$.


 If we define
$\Cal A_i = p_*\Cal O(i)$, then for $i\gg 0$ each $\Cal A_i$ is
locally free of finite rank on $S$, so for $0\ll p \leq q$
the relative graded ideal space $J(k, \Cal A_{[p,q]}/S)$ is
defined. The arguments of this section relativize immediately
to give the following fact.

\proclaim{(1.4.1) Theorem} If $k_i = h^{\Cal O_X}(i)-h(i)$, then
for $0\ll p\ll q$ the natural morphism
$Hilb_h(X/S)\to J(k, \Cal A_{[p,q]}/S)$ is an isomorphism
of $S$-schemes. 
\endproclaim

\vfill\eject

\heading {2. Review of cotangent complexes}\endheading

In this section we review, for future reference, the formalism of
relative cotangent complexes. Our approach, equivalent to
that of Illusie \cite{I}, is based on dg-resolutions,
rather than simplicial resolutions. From now on we assume
that our base field $\Bbb K$ has characteristic 0.

\vskip .3cm

\noindent {\bf (2.1) Generalities.} 
Let $f: X\to Y$ be a morphism of quasiprojective dg-schemes. 
By Theorem 2.7.6 of \cite{CK}, we can factor $f$ into a composition
$X\buildrel i\over\hookrightarrow \widetilde X \buildrel \widetilde f
\over\rightarrow Y$, where $i$ is a quasiisomorphic
closed embedding and $\widetilde f$ is a smooth morphism. 
The relative (derived) cotangent complex $L\Omega^{1\bullet}_{X/Y}$
is, by definition, the dg-sheaf $\Omega^{1\bullet}_{\widetilde X/Y}$
of relative 1-forms for $\widetilde X\to Y$. This is a dg-sheaf
(in fact, a dg-bundle, in the sense of \cite{CK}, Def. 2.3.2)
on $\widetilde X$. 
It is well-defined up to quasiisomorphism
in the sense explained in Proposition 2.7.7 of \cite{CK}. We will denote
$$\Bbb L_{X/Y}^\bullet = i^*\Omega^{1\bullet}_{\widetilde X/Y}.$$
This is a quasi-free dg-sheaf on $X$, well defined up to
quasiisomorphism. Note also that we have a natural projection
(restriction map) $L\Omega^{1\bullet}_{X/Y}\to i_*\Bbb L^\bullet_{X/Y}$
of dg-sheaves on $\widetilde X$. This is a quasiisomorphism, as one
sees from the Eilenberg-Moore spectral sequence (\cite{CK}, Prop.2.4.1). 

\vskip .2cm

\noindent {\bf (2.1.1) Examples.} (a) Note, in particular, the
case when $X,Y$ are usual schemes (trivial dg-structure).
In this case the formalism of relative cotangent
complexes was developed by Illusie \cite{I}. It produces
a simplicial sheaf $L^{\text{Ill}}_\bullet(X/Y)$ of
$\Cal O_X$-modules by using a simplicial  resolution
of $\Cal O_X$ as an $f^{-1}\Cal O_Y$-algebra.
It is not hard to see that the $\Bbb Z_-$-graded cochain
complex associated to  $L^{\text{Ill}}_\bullet(X/Y)$
is quasiisomorphic to $\Bbb L^{\bullet}_{X/Y}$. 
More precisely, we need to use the equivalence between
simplicial commutative algebras and $\Bbb Z_-$-graded
dg-algebras over a field of charcteristic 0. 

\vskip .1cm

(b) If $f: X\to Y$ is a smooth map of usual quasiprojective
schemes, then  we can take $\widetilde X=X$ and
$L\Omega^{1\bullet}_{X/Y} = \Bbb L^\bullet_{X/Y} = \Omega^1_{X/Y}$ is the
usual sheaf of relative 1-forms.

\vskip .1cm

(c) Let $Y$ be a smooth quasiprojective variety
and $f: X\hookrightarrow Y$ be the embedding
of a local complete intersection. If $\Cal I\i \Cal O_Y$
is the sheaf of ideals of $X$, then $\Cal I/\Cal I^2$
is a locally free sheaf on $X$, called the conormal bundle
of $X$ in $Y$ and denoted $\Cal N^*_{X/Y}$. In this
case it is easy to see by using Koszul resolutions
that $\Bbb L^{\bullet}_{X/Y}$ is quasiisomorphic
to $\Cal N^*_{X/Y}$ placed in degree $(-1)$. 

\vskip .1cm

(d) If $\phi: A\to B$ is a morphism of finitely generated
commutative algebras, we have a morphism
of affine schemes $\phi^*: \text{Spec}(B)\to\text{Spec}(A)$.
We will use the notation $L\Omega^{1\bullet}_{B/A}$ 
for the relative cotangent complex of this morphism.
By definition. $L\Omega^{1\bullet}_{B/A} = \Omega^{1\bullet}_{\widetilde B/A}$
where $\widetilde B\to B$ is an $A$-algebra resolution of $B$ such that
$\op{Spec}(\widetilde B)\to \op{Spec}(A)$ is smooth
(in particular, $\widetilde B$ is assumed to have finitely many generators in
each degree). Accordingly, we write
$\Bbb L^\bullet_{B/A} = \Omega^{1\bullet}_{\widetilde B/A}\otimes_{\widetilde B}
B$. This is a complex of $B$-modules well defined up to quasiisomorphism. 
 
In fact, in the affine case we can use more general types of
resolutions, not necessarily finitely
generated. More precisely, an $A$-dg-algebra $F$ is called
quasi-free, if $F^\sharp$ is free over $A^\sharp$ as a graded
commutative algebra.

\proclaim{(2.1.2) Proposition} If, in the situation of Example
2.1.1(d), $F\to B$ is any resolution of $B$ by a quasi-free $A$-dg-algebra,
then $L\Omega^{1\bullet}(B/A)$ is quasiisomorphic to
$\Omega^{1\bullet}(F/A)$, the dg-module of relative 1-forms,
and hence $\Bbb L^\bullet_{B/A} = \Omega^{1\bullet}(F/A)\otimes_F B$. 
\endproclaim

\noindent {\sl Proof:} The  simplicial approach of Illusie, see
Example 2.1.1(a), works without regard to finiteness conditions.
The simplicial module $L_\bullet^{\text{Ill}}(B/A)$ is obtained
by forming a simplicial free resolution, namely the 
bar-resolution, of $B$ as an $A$-algebra.
In order to relate this to our definition, we should
mix the two approaches by defining $L_\bullet^{\text{Ill}}$
for dg-algebras (which is done in an obvious way) and then
prove that for a quasi-free $A$-dg-algebra $F$ the dg-module
obtained from  a simplicial dg-module $L_\bullet^{\text{Ill}}(F/A)$
is quasiisomorphic to $\Bbb L^\bullet_{F/A}$. A spectral sequence
argument reduces this statement to the case when $F$ has trivial
differential, so is just a free graded commutative
algebra over $A$. In this case the proof is similar to that for
a free ungraded algebra which is well known \cite{I}. 

\vskip .3cm

\noindent {\bf (2.2) The graded version.} In this paper
we will allow the dg-algebras (in particular, commutative
algebras with trivial dg-structure) to have an extra
grading indicated in the subscript: $A=\bigoplus_i A_i$,
cf. \S 1 and not contributing to the Koszul sign rule. 
We will call the upper degree cohomological and
lower degree projective.

\vskip .1cm

If $A = \bigoplus A_i$ is a finitely generated
graded commutative algebra (with trivial cohomological
grading), then we denote, as usual, by $\text{Proj}(A)$
the projective spectrum of $A$ and by $\Cal O(1)$
the standard invertible sheaf on it. For a graded $A$-module
$M=\bigoplus M_i$ we denote by $\text{Sh}(M)$ the quasicoherent
sheaf on $\text{Proj}(A)$ obtained from $M$ by
localization. For a quasicoherent sheaf $\Cal F$ on
$\text{Proj}(A)$ we denote by $\text{Mod}(\Cal F)$ the graded
$A$-module $\bigoplus H^0(\text{Proj}(A), \Cal F(i))$. 
It is well known that $\text{Sh}(\text{Mod}(\Cal F))=\Cal F$
while in general there is only a map
$M\to \text{Mod}(\text{Sh}(M))$ which, for $M$ finitely
generated, is an isomorphism in sufficiently
high degrees. 

\vskip .1cm

We now extend this to dg-algebras.  
Let $A$ be  a bigraded dg-algebra as above  with finitely many
generators in each cohomological degree. Then we
have a dg-scheme $\text{Proj}(A) = (\text{Proj}(A^0_\bullet), 
\Cal O^\bullet)$,
where $\Cal O^i$ is the sheaf $\text{Sh}(A^i_\bullet)$,
where $A^i_\bullet$ is regarded as an $A^0_\bullet$-module.
Recall \cite{CK} that a dg-scheme $X$ is called projective
if $X^0$ is projective and each $\Cal O_X^i$ is coherent. From the
above we easily get the following.

\proclaim {(2.2.1) Proposition}
 Any projective dg-scheme $X$ can be obtained as $\roman{Proj}(A)$,
if we take $A^i = \roman{Mod}(\Cal O_X^i)$. 
\endproclaim

For $X=\text{Proj}(A)$ and a bigraded $A$-module $M$
we have a quasicoherent dg-sheaf $\text{Sh}(M)$ on $X$; explicitly,
we take $\text{Sh}(M)^i = \text{Sh}(M^i_\bullet)$. Conversely,
for a quasicoherent dg-sheaf $\Cal F^\bullet$
 on $X$ we form a bigraded $A$-module
$\text{Mod}(\Cal F^\bullet)$ defined by
$\text{Mod}(\Cal F^\bullet)^i_\bullet=\text{Mod}(\Cal F^i)$.
It is clear that $\text{Sh}(\text{Mod}(\Cal F))=\Cal F$.

\vskip .2cm

Let $A=\bigoplus A_i$ and $B=\bigoplus B_i$ be two bigraded
commutative dg-algebras 
and $\phi: A\to B$ be a morphism of dg-algebras
preserving the  grading. In this case the relative
cotangent complex $L\Omega^{1\bullet}_{B/A}$ also
acquires extra grading, induced by those in $A, B$. 
According to our convention, we denote
the $i$th graded component of this complex
by $L\Omega^{1\bullet}_{B/A, i}$. 

\vskip .1cm

Denote $X=\text{Proj}(B)$, $Y=\text{Proj}(A)$ and suppose 
 that we have a morphism $f: X\to Y$
induced by a morphism $\phi: A\to B$ of bigraded dg-algebras. 
We would like to compare $L\Omega^{1\bullet}_{A/B}$ with
$L\Omega^{1\bullet}_{X/Y}$. To construct the former,
it is sufficient to form a free bigraded resolution $\widetilde B$ of $B$
as an $A$-algebra and we can assume that $\widetilde B$ has finitely
many generators in each cohomological degree. Then, we can think
of $L\Omega^{1\bullet}_{A/B} = \Omega^{1\bullet}_{\widetilde B/A}$
 as a bigraded $\widetilde B$-module. Geometrically, the dg-scheme
$\widetilde X = \text{Proj}(\widetilde B)$ is a kind of projective space over $Y=
\text{Proj}(A)$; in particular, we get a factorization of $f$ as
$X=\text{Proj}(B) \i \widetilde X\buildrel {\widetilde f}\over \to Y$,
and this factorization can be used to construct the relative
tangent complex $L\Omega^{1\bullet}_{X/Y}=\Omega^{1\bullet}_{\widetilde X/Y}$
and its restriction $\Bbb L^\bullet_{X/Y}$.

\proclaim{(2.2.2) Theorem} We have the identification of
coherent dg-sheaves on $\widetilde X$:
$$\roman{Sh}(L\Omega^{1\bullet}_{B/A}) = L\Omega^{1\bullet}_{X/Y}.$$
\endproclaim

\noindent {\sl Proof:} First, assume that $X$ and 
$Y$ are ordinary (not dg) smooth
varieties and $f$ is a closed embedding. In this case
$\Bbb L^{\bullet}_{X/Y}$ is equal to 
$\Cal N^*_{X/Y}[1]=\Cal I/\Cal I^2[1]$, where
$\Cal I\subset\Cal O_Y$ is the sheaf of ideals of $X$.
Let us denote by $0$ the ``origin'' in $\text{Spec} (A)$ and 
$\text{Spec}(B)$,
i.e., the point corresponding to the ideal of elements of
positive degree. Then, up to modules supported at $0$, we
have that $L\Omega^{1\bullet}_{B/A}$ is quasiisomorphic to $I/I^2[1]$,
where $I=\text{Ker}(\phi)$. So our statement in the particular case
we consider follows from the next easy lemma.

\proclaim{(2.2.3) Lemma} We have $\roman{Sh}(I^{\nu})=\Cal I^{\nu}$,
for each $\nu\geq 0$.
\endproclaim
 
Now consider the general case and form the fiber product
$\widetilde{X}^2:=\widetilde{X}\times_{Y}\widetilde{X}$, with 
diagonal $\Delta:\widetilde{X}\longrightarrow\widetilde{X}^2$.
Let $\Cal I\subset \Cal O_{\widetilde{X}^2}$ be the dg-ideal
of the diagonal. Then $L\Omega^{1\bullet}_{X/Y}=\Cal I/\Cal I^2$.
Similarly, consider the dg-algebra ${\widetilde B}^2:={\widetilde B}\otimes_A
{\widetilde B}$ and let $I$ be the kernel of the multiplication
homomorphism ${\widetilde B}^2\longrightarrow{\widetilde B}$. Then 
$L\Omega^{1\bullet}_{B/A}$ is quasiisomorphic to $I/I^2$, up to
dg-modules supported at the origin. We conclude the argument
by using an obvious dg-version of Lemma 2.2.3.

\vskip .3cm

\noindent{\bf (2.3) Cotangent complex and Harrison complex.} 
Recall an explicit construction of the absolute cotangent complex
in the affine case via Harrison chains \cite{Lo}. 

Let $A$ be a commutative $\Bbb K$-algebra. We denote by
$\text{FCoLie}(A[1])$ the free graded Lie coalgebra
cogenerated by the graded vector space $A[1]$. Thus
$\text{FCoLie}(A[1]) = A\oplus S^2(A) \oplus ...$
 Let $d$ be the unique differential on $\text{FCoLie}(A[1])$
which is compatible with the Lie coalgebra
structure and on the space of cogenerators is given by the
multiplication $S^2(A)\to A$. The resulting complex
(dg-Lie coalgebra)
$$(\text{FCoLie}(A[1]), d) = \text{Harr}_\bullet(A, \Bbb K)$$
is known as the Harrison chain complex of $A$ with coefficients
in $\Bbb K$. 
More generally, if $M$ is an $A$-module,
then the graded vector space
$$\text{Harr}_\bullet(A, M) = \text{FCoLie}(A[1])\otimes_{\Bbb K} M$$
has a natural differential making it into a dg-comodule over
$\text{Harr}_\bullet(A, \Bbb K)$. It is known as the Harrison chain
complex with coefficients in $M$. Dually, we have the Harrison
cochain complex
$$\text{Harr}^\bullet(A, M) = \text{Hom}_{\Bbb K}(\text{FCoLie}(A[1]), M)=$$
$$\biggl\{ \text{Hom}_{\Bbb K}(A, M)\buildrel \delta\over\lra 
\text{Hom}_{\Bbb K}(S^2(A),M)\lra ...\biggr\}$$
with the first differential acting by the standard formula
$$(\delta f) (a\cdot b) = f(ab) -a f(b) - bf(a).$$
This means that
$$H^1 \text{Harr}^\bullet(A, M) = \text{Der}(A, M), \quad
H^{-1}\text{Harr}(A, M) = \Omega^1_{A/k}\otimes M.$$
Now, the standard property of the Harrison complex  is as follows.

\proclaim{(2.3.1) Theorem} The complex $\roman{Harr}_\bullet(A, A)$
is quasiisomorphic to $\Bbb L^{\bullet}_{A/k}$.
\endproclaim

\noindent {\sl Proof:} This follows from Proposition 2.1.2, applied to
a particular quasi-free resolution of $A$, namely the commutative
 bar-resolution
$$S\Cal L(A) :=
\biggl(\text{FCoLie}(A[1])[-1]\biggr) \buildrel \sim\over\lra A.
\leqno (2.3.2)$$

\heading{3. A finite-dimensional model}\endheading

Let $A$ be a finite dimensional commutative algebra over $\Bbb K$,
possibly without unit. For $k>0$, let $J(k,A)$ be the scheme 
of $k$-dimensional ideals of $A$ (see (1.2)). In this section
we construct a dg-manifold $RJ(k,A)$, with $\pi_0(RJ(k,A))=J(k,A)$
and $$H^iT^{\bullet}_{[I]}RJ(k,A)=\operatorname{Ext}^{i+1}_{A/I}
(\Bbb L^{\bullet}_{(A/I)/A)},\; A/I),$$
for any $\Bbb K$-point $I\in J(k,A)$. As in \cite{CK}, this
is achieved by representing $J(k,A)$ in terms of two abstract
constructions, which we now explain.

\subhead{(3.1) Two constructions}\endsubhead 
(3.1.1) Let $W$ be a finite dimensional vector space over $\Bbb K$. 
We consider the subscheme $CA(W)$ of $\op{Hom}(S^2W,W)$
formed by all commutative, associative multiplications
on $W$. It is clear that for any $\Bbb K$-point  $[\mu]\in CA(W)$,
represented by $\mu:S^2W\lra W$, we have that 
$$T_{[\mu]}CA(W)=Z^2_{\text{Harr}}(W,W).$$
Similarly, if $W$ is a vector bundle over a scheme $S$, we
get the relative space of algebra structures $CA(W)\lra S$.

\vskip .1cm

(3.1.2) Let $S$ be a scheme and let $A$ and $B$ be vector bundles
over $S$ which are made into commutative $\Cal O_S$-algebras. Let
$f:A\lra B$ be morphism of $\Cal O_S$-modules. The {\it homomorphicity
locus} $M_f$ is informally the locus of points $s\in S$ such that
the morphism of fibers $f_s:A_s\lra B_s$ is an algebra homomorphism.
More precisely, $M_f$ is the fiber product

$$\minCDarrowwidth{7 mm}
\CD
 M_f@>>> S\\
@VVV @VV f V\\
 Hom_{\Cal Com}(A,B) @>>> 
\mid Hom_{\Cal O_S}(A,B)\mid 
\endCD
$$
Here $|E|$ means the total space of a vector bundle $E$
and $Hom_{\Cal Com}$ is the subscheme formed by homomorphisms of commutative
algebras. 

\vskip .1cm

We apply these constructions to the following situation. We take $A$ to
be a finite-dimensional algebra, $G(k, A)$ the Grassmannian of $k$-dimensional
subspaces, $\tilde V$ the tautological rank $k$-bundle over $G(k, A)$,
whose fiber over a point $[V]\in G(k, A)$ corresponding to a subspace
$V\subset A$, is $V$. Thus $\tilde V$ is a subbundle in $A$ (the trivial
bundle with fiber $A$). Let $A/\tilde V$ be the quotient bundle. 
Consider the relative scheme of algebra structures on the fibers of this bundle:
$$CA(A/\tilde V)\buildrel q\over \lra G(k, A).$$
On this scheme we have a canonical morphism of vector bundles
$$f: A\to q^*(A/\tilde V).$$
These vector bundles are in fact commutative 
$\Cal O_{CA(A/\tilde V)}$-algebras.

\proclaim {(3.1.3) Proposition} The homomorphicity locus $M_f$
is isomorphic to $J(k, A)$.
\endproclaim

\noindent {\sl Proof:} Almost a tautology. Namely, a $k$-linear
subspace $I\subset A$ is an ideal if and only if $A/I$ is made
into a commutative algebra so that the natural projection
$A\to A/I$ is an algebra homomorphism. 

\vskip .1cm

Now we construct the derived variety of ideals by developing
the derived versions of each of the two above constructions.

\vskip .3cm

\noindent {\bf (3.2) Derived space of algebra structures.} Our approach
is similar to that of C. Rezk \cite{Re} and  parallels
our earlier construction \cite{CK}
for structures of a module, not an algebra. 

First, we generalize the concept of the space of algebra structures.
Let $\Cal P = \{ \Cal P(n), n\geq 0\}$ be a $\Bbb K$-linear operad \cite{GK}.
Thus, by definition, each $\Cal P(n)$ is acted upon by the symmetric
group $S_n$ and $\Cal P$ is equipped with composition maps
$$\Cal P(n)\otimes \Cal P(a_1)\otimes ... \otimes\Cal P(a_n)
\to\Cal P(a_1 +... + a_n)\leqno (3.2.1)$$
satisfying the axioms of \cite {May}. For example, if $W$ is a vector space,
then we have the endomorphism operad $\Cal E_W$
with $\Cal E_W(n) =  \text{Hom}(W^{\otimes n}, W)$ and the
maps (3.2.1) given by composition of multilinear maps.
  A $\Cal P$-algebra structure on $W$ is
a  morphism of operads $\Cal P\to\Cal E_W$. For example, the case
of ordinary commutative algebras (possibly without unit)
corresponds to the case when $\Cal P$ is the commutative
operad $\Cal Com$ with $\Cal Com(n)=\Bbb K$ with trivial $S_n$-action,
for any $n\geq 1$. 

Let now $W$ be a finite-dimensional vector space and $\Cal P$ be
any operad. In this case we have a scheme $\Cal PAlg(W)$ parametrizing
$\Cal P$-algebra structures on $W$. It is realized as a closed
subscheme in the (possibly infinite-dimensional) affine space:
$$\Cal PAlg(W) \subset \prod_n |\text{Hom}_{S_n}(\Cal P(n), \Cal E_W(n))|,
\leqno (3.2.2)$$ 
  given by the equations of ``compatibility with operad structures".
Note that even when $(\Cal P(n))$ is infinite dimensional, 
$|\text{Hom}_{S_n}(\Cal P(n), \Cal E_W(n))|$ still makes sense as a scheme,
the spectrum of a polynomial ring in infinitely many variables. 

There is a case when the scheme of $\Cal P$-algebra structures
can be found explicitly. This is the case of free operads, which
we now recall. Let $E$ be an $S$-module \cite{GK}, i.e., a collection
$E=\{E(n), n\geq 0\}$ of vector spaces together with an $S_n$-action
on $E(n)$ given for each $n$. To $E$, there corresponds
the free operad $F_E$ characterized by the condition that for any
other operad $\Cal P$
$$\text{Hom}_{\text{Operads}}(F_E, \Cal P) = \text{Hom}_{S\text{-modules}}
(E, \Cal P).\leqno (3.2.3)$$
It follows that
$$F_EAlg(W) = \prod_n |\text{Hom}_{S_n}(E(n), \text{Hom}(W^{\otimes n}, W))|
\leqno (3.2.4)$$
is a (possibly infinite-dimensional) affine space. 

We now extend the above formalism to $\Bbb Z_-$-graded
dg-operads, i.e., operads in the symmetric monoidal category
of $\Bbb Z_-$-graded cochain complexes. For such an operad $\Cal P$
by $\Cal P_\sharp$ we will denote the graded operad which is the
same as $\Cal P$ but with differential forgotten.
We say that  a dg-operad $\Cal F$ is quasi-free, if $\Cal F_\sharp$
is free in the above sense. 

\proclaim {(3.2.5) Proposition} Let $W$ be a finite-dimensional $\Bbb K$-vector
space and $\Cal P$ be any $\Bbb Z_-$-graded dg-operad. Then there exists
an affine dg-scheme $\Cal PAlg(W)$ such that for any commutative
dg-algebra $A$ morphisms $\text{Spec}(A)\to \Cal PAlg(W)$ are in bijection
with $A\otimes_{\Bbb K}\Cal P$-algebra structures in the $A$-module
$A\otimes_{\Bbb K}W$. 
\endproclaim

\noindent {\sl Proof:} First, consider the case when $\Cal P$ has trivial
differential (but possibly nontrivial grading). Then, we have
a version of the embedding (3.2.2) with the RHS being the  graded affine
space, namely the spectrum of the $\Bbb Z_-$-graded commutative algebra   
$$\bigotimes_n S^\bullet ((\Cal P(n)\otimes \Cal E_W(n)^*)
^{S_n}).\leqno (3.2.6)$$
The coordinate ring of the
 LHS, i.e., of the sought-for scheme $\Cal PAlg(W)$ is obtained by
quotienting (3.2.6) by the graded ideal formed by the
  equations
expressing compatibility with operad structures. The graded scheme
 $\Cal PAlg(W)$ is, by construction, functorial in $\Cal P$ and $W$. 

Next, suppose that the differential in $\Cal P$ is nontrivial and
consider the graded scheme $\Cal P_\sharp Alg(W)$. The differential
$d$ in $\Cal P$ can be regarded as an odd infinitesimal automorphism of
$\Cal P_\sharp$, see, e.g., \cite{Ka}, Prop. 1.1.1.
 So the graded scheme $\Cal P_\sharp Alg(W)$
inherits this action by naturality, i.e., we get a dg-scheme
which we denote $\Cal PAlg(W)$. The determination of the functor
represented by $\Cal PAlg(W)$ is straightforward.

\proclaim{(3.2.7) Proposition} (a) Any $\Bbb Z_-$-graded dg-operad
$\Cal P$ possesses a quasi-free resolution $\Cal F$. Moreover, if each
graded piece of each $\Cal P(n)$ is finite-dimensional, we can
find $\Cal F$ so that $\Cal F_\sharp = F_E$ and
each graded component of each $E(n)$ is finite-dimensional.

 (b) If $q: \Cal F_1\to\Cal F_2$ is
a quasi-isomorphism of quasi-free dg-operads, then for any
finite-dimensional vector space $W$ the induced morphism of
dg-schemes $q^*: \Cal F_2 Alg(W)\to\Cal F_1Alg(W)$
is a quasiisomorphism. 
\endproclaim

 Part (a) is clear, and part (b) will be proved  in (3.4). Assuming the proposition, we can give the following definition.

\proclaim{(3.2.8) Definition}
Let $\Cal P$ be any $\Bbb Z_-$-graded dg-operad and $W$ a finite-dimensional
vector space. Then the derived space of $\Cal P$-actions
on $W$ is defined to be
$$R\Cal PAlg (W) = \Cal FAlg(W)$$
where $\Cal F\to\Cal P$ is any quasi-free resolution. 
\endproclaim

\vskip .3cm

\noindent {\bf (3.3) The bar-resolution
for operads.}
We now describe a particular functorial quasi-free resolution 
of any dg-operad
\cite{GK}. 

First of all, if $E^\bullet = \{E^\bullet(n)\}$ is any graded $S$-module,
its suspension is the graded $S$-module $\Sigma E^\bullet$
defined by
$$(\Sigma E^\bullet)(n) = E^\bullet(n)[1-n] \otimes\op{sgn}_n,
\leqno (3.3.1)$$
where $\op{sgn}_n$ is the sign representation of $S_n$,
see \cite{GK}. The inverse functor will be denoted by $\Sigma^{-1}$. 
Let now $\Cal P$ be any $\Bbb Z_-$-graded dg-operad. 

\proclaim {(3.3.2) Proposition} (a) For any $\Bbb Z_-$-graded
cooperad $\Cal C$
the free operad $F_{\Sigma^{-1}\Cal C}$ has a natural differential
$d'$ (making it into a dg-operad) which on the cogenerators is induced by
the cocomposition in $\Cal C$.

(b) If $\Cal C$ is a $\Bbb Z_-$-graded dg-cooperad with differential $d_{\Cal C}$,
then the induced differential $d''$ on $F_{\Sigma^{-1}\Cal C_\sharp}$
commutes with $d'$ and the total differential $d=d'+d''$
makes $F_{\Sigma^{-1}\Cal C}$ into a dg-operad
which we denote $\op{Cobar}(\Cal C)$. 
\endproclaim

This construction is known as the cobar-construction of the
cooperad $\Cal C$, see \cite{GK} \cite{GJ}. Similarly, we have the bar-construction
for operads.

\proclaim{(3.3.3) Proposition} (a) For any $\Bbb Z_-$-graded operad
$\Cal P$ the free cooperad $\op{Bar}(\Cal P)$ generated by $\Sigma\Cal P$
has a natural differential $d'$ (making it into a dg-cooperad)
which on the space of generators is induced by the composition in
$\Cal P$.

(b) If $\Cal P$ is a $\Bbb Z_-$-graded dg-operad with differential
$d_{\Cal P}$, then the induced differential $d''$
on $\op{Bar}(\Cal P_\sharp)$ commutes with $d'$ and the
total differential $d=d'+d''$ makes it into a dg-cooperad which
we denote $\op{Bar}(\Cal P)$.

\endproclaim

 By construction, Bar and Cobar are functors from dg-operads
to dg-cooperads and back, and it is easy to see
that these functors take quasi-isomorphisms to quasi-isomorphisms. 

We now define the {\it bar-resolution} of a dg-operad $\Cal P$ to be
$$\Cal B(\Cal P) = \op{Cobar}(\op{Bar}(\Cal P)).\leqno (3.3.4)$$
Thus,  $\Cal B(\Cal P)$ is quasi-free and functorial 
in $\Cal P$.

\proclaim{(3.3.5) Proposition} There is a natural quasiisomorphism
$\alpha: \Cal B(\Cal P)\to\Cal P$. 
\endproclaim

For a vector space $W$ this particular quasi-free resolution
gives a particular model for the derived space
of $\Cal P$-algebra structures which we denote
$$\widetilde R\Cal PAlg(W) = \Cal B(\Cal P)Alg(W).\leqno (3.3.6)$$
A $\Cal B(\Cal P)$-algebra $W$ is sometimes called
a homotopy $\Cal P$-algebra.  More precisely, we have
maps $\Cal P(n)\otimes W^{\otimes n}\to W$ which satisfy the
axioms of an operad action only up to higher homotopies. 

In the same way as Proposition 3.5.3 of \cite{CK}, we prove the following.

\proclaim {(3.3.7) Proposition} For any $\Bbb Z_-$-graded dg-operad
$\Cal P$ we have natural convergent spectral sequences
$$E_1 = H^\bullet \Bbb K [\widetilde R\Cal P_\sharp Alg(W)] \quad \Rightarrow
\quad H^\bullet \Bbb K[\widetilde R\Cal P Alg(W)]; \leqno \op{(a)}$$
$$E_2 = \Bbb K[\widetilde R H^\bullet(\Cal P) Alg(W)]
 \quad \Rightarrow
\quad H^\bullet \Bbb K[\widetilde R\Cal P Alg(W)]. \leqno \op{(b)}$$
\endproclaim

\vskip .3cm

\noindent {\bf (3.4) M-homotopies and the proof of Proposition
3.2.7.}
Let $f: \Cal P\to\Cal Q$ be a morphism of dg-operads. Then $f$ makes
$\Cal Q$ into a bimodule over $\Cal P$ in the sense of \cite{Mar}. 
We can, therefore, speak about derivations $D: \Cal P\to \Cal Q$
with respect to this bimodule structure which are collections
of morphisms of complexes $D_n: \Cal P(n)\to\Cal Q(n)$ such that for any
$n, a_1, ..., a_n\in \Bbb Z_+$ and any $p\in \Cal P(n)$,
$ p_i\in \Cal P(a_i)$ we have
$$D_{a_1+...+a_n} \bigl( p(p_1, ..., p_n)\bigr) = (D_n(p))(f(p_1), ..., f(p_n))
+$$
$$+\sum_{i=1}^n f(p)\bigl( f(p_1), ..., D_{a_i}(p_i), ..., f(p_n)\bigr).$$
Here $p(p_1, ..., p_n)$ is the image of $p\otimes p_1\otimes ... \otimes
p_n$ in $\Cal P(a_1 + ... + a_n)$ under the composition map in $\Cal P$.
More generally, we define derivations of degree $d$ by allowing
the $D_n$ to be morphisms of complexes of degree $d$ and
introducing obvious sign factors.  

It is clear that whenever $(f_t: \Cal P\to\Cal Q)_{t\in [0,1]}$
is a smooth family of morphisms, then for each $t$ the
derivative $f'_t = {d\over dt} f_t$ is a derivation $\Cal P\to\Cal Q$
with respect to the bimodule structure given by $f_t$.  

\proclaim{(3.4.1) Definition} An $M$-homotopy between two morphisms
$f, g: \Cal P\to\Cal Q$ is a pair $(f_t, s_t)$, where
$(f_t)_{t\in [0,1]}$ is a smooth family of morphisms of dg-operads
such that $f_0=f$, $f_1=g$ and $(s_t)$ is a smooth family of
degree (-1) derivations $\Cal P\to \Cal Q$ (with respect to the
bimodule structure given by $f_t$) such that $f'_t = [d, s_t]$.
\endproclaim

As in the case for algebras \cite{CK}, we see that two $M$-homotopic
morphisms $\Cal P\to\Cal Q$ induce the same morphism
$H^\bullet(\Cal P)\to H^\bullet(\Cal Q)$. 

\proclaim{(3.4.2) Lemma} Let $\Cal P, \Cal Q$ be $\Bbb Z_-$-graded
dg-operads such that $\Cal P$ is quasifree and each $\Cal Q(n)$
is acyclic in degrees $<0$. Suppose that $f_1, f_2: \Cal P\to\Cal Q$
are two morphisms of dg-operads which, for each $n$, induce
the same morphisms $H^0\Cal P(n)\to H^0\Cal Q(n)$.
Then there exists an $M$-homotopy connecting $f_0$ and $f_1$.
\endproclaim

\noindent{\sl Proof:} This is similar to the proof
of Proposition 3.6.4 of \cite{CK} which deals with
an analogous statement but for associative algebras
instead of operads. More precisely, the
inductive (in the homological degree) procedure
from that proof adapts without difficulty to the
case of operads.

\proclaim {(3.4.3) Corollary} Let $\Cal R = F_E$ be the
free graded operad on a $\Bbb Z_-$-graded $S$-module $E$
(no differential). 
Then the quasiisomorphism $\alpha: \Cal B(\Cal R)\to\Cal R$
from Proposition 3.3.5 has a natural left inverse
$\beta: \Cal R\to\Cal B(\Cal R)$ and the composition
$\beta\alpha$ is $M$-homotopic to the identity of
$\Cal B(\Cal R)$.
\endproclaim

\noindent {\sl Proof:} As $\Cal R$ is free, we can define
a morphism from $\Cal R$ by prescribing its restriction
to  the $S$-module of generators. We define $\beta$ on
$E$ to be the canonical identification of it with
the natural copy of $E$ inside $\Cal R\i F_{\Cal R}\i
\Cal B(\Cal R)$. Now, if $E$ is concentrated in degree 0
(i.e., each $E(n)$ is), then we can apply Lemma 3.4.2 to
$\Cal P=\Cal Q = \Cal B(\Cal R)$. If
$E$ is not concentrated in degree 0, then we notice
that $\Cal B(\Cal P)$ comes from a $\Bbb Z_-\times \Bbb Z_-$-graded
dg-operad with differential of degree (1,0) and then
use the same proof as for Lemma 3.4.2, but with induction
in the second component of the bidegree. 

\vskip .1cm

We now pass to the proof of Proposition 3.2.7. Again, we use
the same method as was used for the proof of Proposition
3.3.6 of \cite{CK} so we will
just outline the main steps. 

First, we prove that if
$\Cal P=F_E$ is free with trivial differential, then
$\widetilde R\Cal PAlg(W)$ is quasiisomorphic to
$\Cal PAlg(W)$. This is done by using Corollary 3.4.3 to $\Cal R=\Cal P$
and noticing that the identity $\alpha\beta = \op{Id}$ and
the M-homotopy $\beta\alpha\sim\op{Id}$ between morphsisms
of operads are inherited in functorial constructions such as 
passing from an operad $\Cal Q$ to the coordinate
algebra $\Bbb K[ \Cal QAlg(W)]$.

Then, to finish the argument,  we prove that for any quasi-free
resolution $p: \Cal F\to\Cal P$ the dg-algebra
$\Bbb K[\Cal FAlg(W)]$ is naturally quasiisomorphic to
$\Bbb K[\widetilde R \Cal PAlg(W)]$
(so is independent, up to quasiisomorphism, on $\Cal F$).
 Indeed, by the above
$\Bbb K[\Cal F_\sharp Alg(W)]$ is quasiisomorphic
to $\Bbb K[\widetilde R \Cal F_\sharp Alg(W)]$.
Then, the first spectral sequence of 3.3.7 implies that 
$\Bbb K[\Cal FAlg(W)]
\to \Bbb K[\widetilde R \Cal FAlg(W)]$ is a quasiisomorphism
while the second sequence implies that 
$\Bbb K[\widetilde R \Cal F Alg(W)] \to 
\Bbb K[\widetilde R \Cal P Alg(W)]$ is a quasiisomorphism.

\vskip .3cm

\noindent {\bf (3.4) The small bar-resolution for $\Cal P = \Cal Com$.}
We now specialize the above discussion to $\Cal P = \Cal Com$
and write $RCA(W) = R\Cal Com Alg (W)$. In addition to the
bar-resolution $\Cal B (\Cal Com)$ we will use another
quasi-free resolution $\Lambda\to\Cal Com$ which we call
the {\it small bar-resolution} and which is defined as follows.

Let $\Cal Lie$ be the operad describing Lie algebras. Denote
by $\Cal Lie^*$ the cooperad formed by the dual spaces
$\Cal Lie(n)^*$. We set $\Lambda = \op{Cobar}(\Cal Lie^*)$. 
Because $\Cal Lie^*$ is quasiisomorphic to $\op{Bar}(\Cal Com)$
(Koszul duality, see \cite{GK}), we find that $\Lambda$ is
quasiisomorphic to $\Cal B(\Cal Com)$; in particular,
it is a quasi-free resolution of $\Cal Com$.

We will denote $\bar RCA(W)$ the particular model for
$RCA(W)$ obtained by using the resolution $\Lambda$.

\proclaim{(3.4.1) Proposition} The affine dg-scheme $\bar RCA(W)$
is a dg-manifold. In fact, the dg-algebra
$\Bbb K[\bar RCA(W)]_\sharp$ is free with finitely many
generators in each degree.
\endproclaim

\noindent {\sl Proof:} For a vector space $W$ a $\Lambda$-algebra structure
on $W$ is the same as a differential $D$ in the free Lie
coalgebra on $W[1]$ satisfying $D^2=0$ and compatible with
the coalgebra structure. So $\Lambda$-algebras are the
same as homotopy commutative algebras in the sense of
Stasheff \cite{St}. Now, a $D$ as before
is defined by the projection of its image on $W[1]$,
and this projection can be arbitrary. So $D$ is defined
by a collection of linear maps $D_i: \op{FCoLie}_i(W[1])\to W$,
$i\geq 2$, $\op{deg}(D_i) = 2-i$, where  $\op{FCoLie}_i(W[1])$ is the
$i$th graded component of $\op{FCoLie}(W[1])$.
Notice that it is finite-dimensional.
Now, the graded algebra $\Bbb K[\bar RCA(W)]_\sharp$
is, by construction, freely generated by
by the  matrix elements of indeterminate maps $D_i$. 

\proclaim{(3.4.2) Proposition} Let $W$ be a finite-dimensional
vector space. Then:

(a) $\pi_0 RCA(W) = CA(W)$.

(b) If $\mu: S^2 W\to W$ represents a $\Bbb K$-point $[\mu]\in CA(W)$,
then
$$H^i T^\bullet_{[\mu]} RCA(W) = \cases Z^2_{\op{Harr}}(W, W), &  i=0;\\
H^{i+2}_{\op{Harr}}(W,W), &  i>0.\endcases $$
\endproclaim
 
\noindent {\sl Proof:} Part (a) is obvious from any construction of $RCA$.
As to (b), the definition of the Harrison complex uses the free
Lie coalgebra on $W[1]$, i.e., exactly the structure involved
in the resolution $\Lambda$.

\vskip .3cm

\noindent {\bf (3.5) The derived homomorphicity locus.} 
Let $A, B$ be two commutative $\Bbb K$-algebras and $\dim(B)< \infty$.
We have then the scheme of homomorphisms
$$\op{Hom}_{\Cal Com}(A, B) \i |\op{Hom}_{\Bbb K}(A, B)|$$
whose points are morphisms of commutative algebras $A\to B$. If $[f]$
is a $\Bbb K$-point represented by a homomorphism $f: A\to B$, then, clearly,
$$T_{[f]} \op{Hom}_{\Cal Com}(A, B) = \op{Der} (A, B), \leqno (3.5.1)$$
where $B$ is made into an $A$-algebra via $f$. 

Let now $A$ be a $\Bbb Z_-$-graded dg-algebra. Then
$|\op{Hom}_{\Bbb K}(A, B)|$ is the $\Bbb Z_-$-graded dg-scheme corresponding to
the $\Bbb Z_+$-graded complex $\op{Hom}_{\Bbb K}(A, B)$ see \cite{CK}, 
n. (2.2.5). The ideal of the graded subscheme
$\op{Hom}_{\Cal Com}(A_\sharp , B) \i |\op{Hom}_{\Bbb K}(A_\sharp, B)|$,
is a differential ideal, so we get the dg-scheme
$\op{Hom}_{\Cal Com}(A, B)$ which parametrizes, in the obvious
sense, dg-algebra homomorphisms. We now define the derived
scheme of homomorphisms to be
$$\op{RHom}_{\Cal Com}(A, B) = \op{Hom}_{\Cal Com}(F, B),\leqno (3.5.2)$$
where $F\to A$ is any quasi-free commutative dg-algebra
resolution. 

\proclaim{(3.5.3) Lemma} The definition of 
$\op{RHom}_{\Cal Com}(A, B)$ is independent, up to quasiisomorphism,
on the choice of $F$.
\endproclaim

\noindent {\sl Proof:} This is done similarly to the proof
of Proposition 3.2.7 and of Proposition 3.3.6 of \cite{CK}, by using the
following ingredients: the concept of $M$-homotopies
for commutative dg-algebras (parallel to the associative
case treated in \cite{CK}) and the canonical quasi-free
resolution $\Cal L(A)= S\bigl(\op{FCoLie}(A[1])[-1]\bigr)\to A$,
which was already used (for ungraded algebras) in (2.3.1). 

\proclaim {(3.5.4) Proposition} Let $A, B$ be ungraded
commutative algebras, $\dim(B)<\infty$. Then:

(a) $\pi_0 \op{RHom}_{\Cal Com}(A,B) = \op{Hom}_{\Cal Com}(A, B)$.

(b)if $f: A\to B$ is an algebra homomorphism, then
$$H^i T^\bullet_{[f]} \op{RHom}_{\Cal Com}(A,B) = H^{i+1}_{\op{Harr}}
(A, B).$$
\endproclaim

Note that the Harrison cohomology are the higher derived functors
of the space of derivations which, by (3.5.1), is the
tangent
space to $\op{Hom}_{\Cal Com}(A, B)$. 

\vskip .1cm

\noindent {\sl Proof:} (a) is clear and (b) is obtained
by using  the resolution $\Cal L(A)$. 

\vskip .1cm

We now globalize this construction as follows. Let $S$ be a 
dg-scheme and $A, B$ be quasicoherent sheaves of $\Cal O_S$-dg-algebras.
We assume that as sheaves of graded $\Cal O_{S \sharp}$-modules,
$A_\sharp$ and $B_\sharp$ are locally free and, moreover,
generators of $A_\sharp$ are in degrees $\leq 0$ and generators
of $B_\sharp$ are finite in number and have degree 0. 
Then we can apply the construction of the derived
homomorphicity locus relative over $S$, getting a
dg-scheme $R\Cal Hom_{\Cal Com \otimes\Cal O_S}(A, B)$.
It comes equipped with a canonical $S$-map $p$ into the total space
$|\Cal Hom_{\Cal O_S} (A, B)|$. 

Let now $f: A\to B$ be a morphism of $\Cal O_S$-dg-modules.
We define the derived homomorphicity locus $RM_f$
to be   the  derived fiber product (in the sense of n. (2.8) of \cite{CK})
$$\minCDarrowwidth{7 mm}
\CD
 M_f@>>> S\\
@VVV @VV f V\\
 R\Cal Hom_{\Cal Com \otimes\Cal O_S}(A,B) @>>> 
\mid \Cal Hom_{\Cal O_S}(A,B)\mid 
\endCD
$$

Notice that the construction of the derived homomorphicity locus
can be applied to the more general situation when $A$ and
$B$ are sheaves of $\Cal P\otimes \Cal O_S$-dg-algebras,
where $\Cal P$ is any $\Bbb Z_-$-graded dg-operad (instead of $\Cal Com$).
In particular, we will use this construction for $\Cal P=\Lambda$. 

\vskip .3cm

\noindent {\bf (3.6) The derived ideal scheme.}
Let $A$ be a finite-dimensional ungraded commutative $\Bbb K$-algebra.
We now apply the derived versions of the two constructions of (3.1),
from which we use the notation. We first consider the
relative derived space of algebra structures
$$RCA(A/\tilde V)\buildrel q\over\lra G(k, A).\leqno (3.6.1)$$
As explained, this is done using a quasi-free resolution
$\Cal F\to\Cal Com$. By construction, $q^*(A/\tilde V)$ is a
sheaf of $\Cal O_{RCA(A/\tilde V)}\otimes \Cal F$-dg-algebras.
 We have
an $\Cal O_{\bar RCA(A/\tilde V)}$-linear morphism
$$f: A\otimes  \Cal O_{\bar RCA(A/\tilde V)}\to q^*(A/\tilde V)$$
whose source and target are sheaves of $\Cal O_{RCA(A/\tilde V)}\otimes 
\Cal F$-dg-algebras. 

\proclaim {(3.6.2) Definition}
(a) We define 
the derived scheme of ideals $RJ(k, A)$ to be
the derived linearity locus $RM_f$. 

(b) By  $\bar RJ(k, A)$ we denote the particular model
for $RJ(k, A)$ obtained by using the resolution
$\Cal F=\Lambda$ for the operad $\Cal Com$
and the resolution $\Cal L(A)$ for the algebra $A$.

\endproclaim

It is clear from the above that  $RJ(k, A)$ is independent,
up to quasi-isomorphism, on the choice of $\Cal F$,
so $\bar RJ(k, A)$ is  particular representative in
the quasi-isomorphism class of $RJ(k, A)$. By (3.4.1), 
$\bar RJ(k, A)$ is a dg-manifold.

 \proclaim{(3.6.3) Theorem} (a) We have $\pi_0 (RJ(k, A)) = 
J(k, A)$.

(b) If $I\i A$ is an ideal representing a $\Bbb K$-point $[I]\in J(k, A)$,
then
$$H^iT^{\bullet}_{[I]}RJ(k,A)=\operatorname{Ext}^{i+1}_{A/I}
(\Bbb L^{\bullet}_{(A/I)/A},\; A/I).$$
\endproclaim

\noindent {\sl Proof:} Part (a) is clear from the corresponding
statements about the two constructions in (3.1). To prove (b),
recall the transitivity triangle of relative cotangent complexes
 (see, e.g., \cite{Q}, Th. 5.1) corresponding to the morphisms of
rings $\Bbb K\to A\to A/I$: 
$$(A/I)\otimes_A \Bbb L^{\bullet}_{A/\Bbb K}\lra
\Bbb L^{\bullet}_{(A/I)/\Bbb K}\lra \Bbb L^{\bullet}_{(A/I)/A}\lra
(A/I)\otimes_A \Bbb L^{\bullet}_{A/\Bbb K}[1]. \leqno (3.6.4)$$
This is a distinguished triangle in the derived category
of complexes of $A/I$-modules. By applying the functor
$\op{RHom}_{A/I}(-, A/I)$, we get the triangle
$$\op{RHom}_{A/I}(\Bbb L^{\bullet}_{(A/I)/A}, A/I)\to
\op{RHom}_{A/I}( \Bbb L^{\bullet}_{(A/I)/\Bbb K}, A/I)\to
\op{RHom}_{A/I}(\Bbb L^{\bullet}_{A/\Bbb K}, A/I), \leqno (3.6.5)$$
of which the first term has the cohomology appearing in
the statement of (b). 
Notice that the other two terms are directly related with
the derived versions of the two constructions of (3.1). 
Namely, Theorem 2.3.1 together with Proposition 3.4.2(b)
imply that the middle term is quasi-isomorphic to the tangent
space of $RCA(A/\tilde V)$ at the point represented by $I$
and by the canonical algebra structure on $A/I$.
Similarly, Proposition 3.5.4(b) implies that
the right term is quasi-isomorphic to the
tangent space to $RM_f$ at the  same point as above. 
Notice further that the successive application of the two constructions
($RCA$ and $M_f$) gives rise, at the level of tangent spaces,
to the cone of a natural morphism of the tangent
spaces to the dg-schemes given by each construction separately.
Theorem is proved. 

\vskip .2cm

\noindent {\bf (3.6.6) The graded case.}
Suppose that $A=\bigoplus A_n$ has an extra grading, as in (1.2).
Then all the previous constructions are modified straightforwardly
to give the {\it derived scheme of graded ideals} $RJ(k, A)$,
where $k=(k_i)$ is sequence of integers.
As before, we use the notation $\bar RJ(k, A)$
for the particular model obtained by using the
resolutions mentioned in Definition 3.6.2.

 In particular,
for any  graded ideal $I\i A$ representing a point $[I]\in J(k, A)$
 we have
$$H^iT^\bullet_I RJ(k, A) = \op{Ext}^{i+1, 0}_{A/I} (\Bbb L^\bullet_{(A/I)/A},
A/I), \leqno (3.6.7)$$
where $\op{Ext}^{i+1, 0}$ is the degree 0 component of $\op{Ext}^{i+1}$. 

\vskip .3cm

\noindent {\bf (3.7) The relative version.} Let $S$ be a scheme
over $\Bbb K$ and $\Cal A$ be a sheaf of commutative $\Cal O_S$-algebras
(possibly without unit) such that as a sheaf of $\Cal O_S$-modules,
$\Cal A$ is locally free of finite rank. Then the
construction of $\bar RJ(k, A)$ relativizes immediately, giving the relative derived ideal
scheme $\bar RJ(k, \Cal A/S)$ which is a dg-scheme equipped with
a smooth morphism to $S$. To construct it, we first form
the relative Grassmannian $G(k, \Cal A/S)\to S$,
which is equipped with the tautological quotient bundle
$\Cal A/\widetilde V$, and then form $RCA(\Cal A/\widetilde V)$, the derived
space of commutative algebra structures in the fibers of 
$\Cal A/\widetilde V$ and the derived homomorphicity locus
of the natural morphism of sheaves of (homotopy) commutative
algebras on it, as in (3.6). 
By construction (and by the naturality  and the tensor nature of
the resolutions $\Lambda$ for the operad $\Cal Com$ and $\Cal L(\Cal A)$
for the $\Cal O_S$-Algebra $\Cal A$),
we have the following compatibility statement.

\proclaim{(3.7.1) Proposition} (a) We have
$\pi_0 RJ(k, \Cal A/S) = J(k, \Cal A/S)$, see (1.4). 

(b)
Let $\phi: S'\to S$ be a morphism of schemes and  $\Cal A$ a sheaf
of commutative $\Cal O_S$-algebras, which is locally free as a sheaf
of $\Cal O_S$-modules. Then $\bar RJ(k, (\phi^*\Cal A)/S')$ is
isomorphic to the  fiber product of
$\bar RJ(k, \Cal A/S)$ and $S'$ over $S$.
\endproclaim

See \cite{CK}, Sect. (2.8) for the definition of the fiber products
in the dg-category. Note that we consider here the
``straight" fiber product (which in our situation is quasiisomorphic
to the derived one). 

Further, if $\Cal A$ is equipped with an extra grading $\Cal A = \bigoplus
\Cal A_n$, as in (1.2) and (3.6.6) (so that each $\Cal A_n$ is
an $\Cal O_S$-submodule in $\Cal A$)
 and $k= (k_i)$ is a sequence of integers, then we have
the relative derived scheme of graded ideals $\bar RJ(k, \Cal A/S)$
defined in an obvious way. 

\vfill\eject

\heading{4. The derived Hilbert scheme}\endheading

\noindent {\bf (4.1) Graded ideals in truncations of the coordinate
ring.} 
We return to the situation of (1.3), so $A = \bigoplus_i A_i$
is the coordinate algebra of a projective variety $X\i \Bbb P^n$.

We consider the finite-dimensional graded algebra $A_{[p,q]}$ and
the corresponding derived scheme of graded ideals $RJ(k, A_{[p,q]})$.
By Theorem 1.3.1, for $0\ll p\ll q$ we have
$\pi_0 RJ(k, A_{[p,q]}) = Hilb_h(X)$. We would like to show
that the dg-structure on $RJ(k, A_{[p,q]})$ is ``asymptotically correct".
More precisely, let $Z\i X$ be a closed subscheme with Hilbert
polynomial $h$, let $\Cal I\i \Cal O_X$ its
sheaf of ideals and 
  $I = \op{Mod}(\Cal I)$ the corresponding graded ideal in $A$. 
We will prove the following fact.

\proclaim{(4.1.1) Theorem} For  $i>0$, for $p$ large enough
(depending on $i$ and $q$ large enough (depending on $i, p$) and
  every $Z\in Hilb_h(X)$  as before we have
$$ \op{Ext}^i_{\Cal O_Z} (\Bbb L^\bullet_{Z/X}, \Cal O_Z) = 
\op{Ext}^{i, 0}_{A_{[p,q]}/I_{[p,q]}}
(\Bbb L^\bullet_{(A_{[p,q]}/I_{[p,q]})/A_{
[p,q]}}, A_{[p,q]}/I_{[p,q]}).$$
When $A$ varies
in a flat family over a scheme $S$ of finite type, then numbers
$p, q$ can be chosen so as the above conditions hold uniformly
for all $\Bbb K$-points of $S$. 
\endproclaim

\proclaim{(4.1.2) Lemma} For any $i>0$,  for 
$p$ large enough (depending on $i$), 
and any $Z\in Hilb_h(X)$ we have
$$\op{Ext}^i_{\Cal O_Z} (\Bbb L^\bullet_{Z/X}, \Cal O_Z) =
\op{Ext}^{i, 0}_{A_{\geq p}/I_{\geq p}} (\Bbb L^\bullet_{(A_{\geq p}/I_{\geq p})
/A_{\geq p}}, A_{\geq p}/I_{\geq p}).$$
\endproclaim

\noindent {\sl Proof:} First, consider two coherent sheaves $\Cal F, \Cal G$
on $Z$ and let $M, N$ be the corresponding graded $A/I$-modules. 
Then Serre's theorem \cite {Se} implies that 
$$\op{Ext}^i_{\Cal O_Z}(\Cal F, \Cal G) = \ind_p \op{Ext}^{i, 0}_{A/I}
(M_{\geq p},
N_{\geq p}),$$
that this limit is achieved and, moreover, achieved uniformly
if $Z, \Cal F, \Cal G$ run in a flat family over  a scheme of finite type. 
We claim that, further, for $p\gg 0$ (with the same uniformity conditions),
$$\op{Ext}^{i, 0}_{A/I}(M_{\geq p},
N_{\geq p}) = \op{Ext}^{i, 0}_{A_{\geq p}/I_{\geq p}}(M_{\geq p},
N_{\geq p}).$$
Indeed, in the case $i=0$ (when we are dealing with Hom),
this is obtained by using the trick with the Veronese subalgebra
from the proof of Theorem 1.3.2. The case $i>0$ is formally
deduced from this by applying the derived functor.

We now specialize to the case when $\Cal G= \Cal O_Z$
and $\Cal F$ is one of the first $i$ cohomology sheaves
of the complex $\Bbb L^\bullet_{Z/X}$. Since, by Theorem 2.2.2,
$\op{Sh}(\Bbb L^\bullet_{(A/I)/A}) = \Bbb L^\bullet_{Z/X}$,
the lemma follows by a spectral sequence argument. 

\vskip .2cm

In virtue of the lemma, Theorem 4.1.1 would be implied by
the following general proposition applied to the canonical morphism
$A_{\geq p} \to A_{\geq p}/I_{\geq p}$.

\proclaim{(4.1.3) Proposition} Let $f: A\to B$ be a morphism of
finitely generated commutative graded algebras, situated in
projective degrees $\geq 0$ (and cohomological degree 0). Then
for each $i$ and large enough $q$ (depending on $i$)
we have 
 $$\op{Ext}^{i,0}_B (\Bbb L^\bullet_{B/A}, B) = \op{Ext}_{B_{\leq q}}^{i,0}(
L^\bullet_{B_{\leq q}/A_{\leq q}}, B_{\leq q}).$$
Further, for fixed $i$ one can choose $q$ uniformly, if
$A, B, f$ vary in a flat family over a scheme of finite type.
\endproclaim

\noindent {\sl Proof:}
For a bigraded dg-algebra $C$ and a dg-module $M$ over $C$ we write
$$\op{RDer}(C, M) = \op{RHom}_C(\Bbb L^\bullet_{C/\Bbb K}, M).$$
Note that if $C$ has trivial cohomological grading, then
$$H^0 \op{RDer}(C, M) = \op{Der}(C, M)$$
is the space of derivations $C\to M$. We will use the
notation $\op{Der}(C, M)$ to signify the space of derivations
for arbitrary $C, M$ as well. 

Returning to the situation of $f: A\to B$, we have the transitivity
triangle similar to (3.6.4), which we write in the form
$$\op{RHom}_B(\Bbb L^\bullet_{B/A}, B) \to \op{RDer}(B, B) \to
\op{RDer}(A, B).$$
We have a similar triangle for the truncated morphisms $A_{\leq q}
\to B_{\leq q}$. Thus our statement  follows from the next proposition.

\proclaim{(4.1.4) Proposition} Let $M$ be a finitely generated $B$-module
(situated in cohomological degree 0). Then for each $i$ and large
enough $q$ (depending on $i$) we have
$$H^{i, 0} \op{RDer}(B, M) = H^{i, 0} \op{RDer}(B_{\leq q}, M_{\leq q}),$$
and $q$ can be choosen uniformly for a flat family of $B, M$ over a scheme
of finite type. 
\endproclaim

\noindent {\sl Proof:} Let $F\to B$ be a quasi-free resolution. Thus $F$
is a bigraded dg-algebra situated in cohomological degrees $\leq 0$ 
and we can assume that $F$ has finitely many generators in each
cohomological degree. By definition,
$$\op{RDer}(B, M) = \op{Der}(F, M).$$
We first prove the following fact.

\proclaim{(4.1.5) Lemma} For each $i$ and for $q\gg 0$ (depending on $i$)
we have
$$H^{i, 0} \op{Der}(F, M) = H^{i, 0} \op{Der}(F_{\leq q}, M_{\leq q}).$$
\endproclaim

\noindent {\sl Proof:} A derivation $D: F\to M$ is uniquely
determined by its values on the generators of $F$, which
can be chosen in an arbitrary way. So if we choose $q$ to be
greater than the projective degrees of the generators
of $F$ in cohomological degrees $0, -1, ..., -i-1$, then the
complex $\op{Der}^0(F, M)$ of derivations preserving the\
projective degrees, coincides, in degrees $\leq i+1$, with
the similar complex $\op{Der}^0(F_{\leq q}, M_{\leq q})$
and hence the two complexes have the same cohomology in degrees $0, 1, ..., i$.
The uniformity in this situation is clear.

\proclaim {(4.1.6) Lemma} 
For each $i, q$ we have
$$H^{i, 0} \op{Der}(F_{\leq q}, M_{\leq q}) = H^{i, 0} \op{RDer}(F_{\leq q},
M_{\leq q}).$$
\endproclaim

\noindent {\sl Proof:} To find $\op{RDer}(F_{\leq q}, M_{\leq q})$,
we need a quasi-free resolution $G$ of $F_{\leq q}$.
To construct such a $G$, we start from $F$ and then add new generators
to $F$ so as, first, to make all the cocycles of $F$ in projective
degrees $>q$ into coboundaries and then, inductively, to similarly
kill any new cocycles that arise. In this way we get a factorization
$$F\hookrightarrow G\buildrel \op{qis}\over \lra F_{\leq q},$$
and any new generator which is added in forming $G$ starting from $F$
is in projective degree $>q$. This means that
$$\op{RDer}(F_{\leq q}, M_{\leq q}) := \op{Der}(G, M_{\leq q}) = 
\op{Der}(F, M_{\leq q}) = \op{Der}(F_{\leq q}, M_{\leq q}),$$
because any of the new generators must necessarily go to 0
in $M_{\leq q}$ by  (projective) degree reasons. This proves the lemma.
Proposition 4.1.4 follows. 

This finishes the proof of Theorem 4.1.1.

\vskip .3cm

\noindent {\bf (4.2) The approximations to $RHilb$.} We continue to
work in the situation of (4.1) and would like to use
the dg-manifolds $RJ(k, A_{[p,q]})$ as approximations
to the derived Hilbert scheme. 

Let $X=(X^0, \Cal O^\bullet_X)$ be a dg-scheme, so that
we have the sheaf of graded algebras $\underline{H}^\bullet(\Cal O^\bullet_X)$
on $\pi_0(X) = \op{Spec}\underline{H}^0(\Cal O_X)$. Let
$\underline{H}^{\geq -m}(\Cal O^\bullet_X)$, $m\geq 0$, denote
the truncation of this sheaf obtained by disregarding
the cohomology in degrees $<-m$. We denote
by $X_h^{\leq m}$ the graded scheme 
$(\pi_0(X), \underline{H}^{\geq -m}(\Cal O^\bullet_X))$.
We call a morphism $f: X\to Y$ of dg-schemes
 an $m$-quasiisomorphism if the induced morphism of graded schemes
$f_h^{\leq m}: X_h^{\leq m}\to Y_h^{\leq m}$ is an isomorphism.

\proclaim {(4.2.1) Proposition} Let $f: M\to N$ be a morphism
of dg-manifolds. Then the following are equivalent:

(i) $f$ is an $m$-quasiisomorphism;

(ii) $\pi_0(f): \pi_0(M)\to\pi_0(N)$ is an isomorphism and
for any field extension $\Bbb F\supset \Bbb K$ and 
any $\Bbb F$-point $x\in M$ the differential $d_xf$
indces isomorphisms $\pi_i(M, x) \to\pi_i(N, f(x))$ for $i\geq -m$.

\endproclaim

Recall that $\pi_i(M, x) = H^{-i}(T^\bullet_x M)$. The proposition
is proved in the same way as Proposition 2.5.9 of \cite{CK}. 

\vskip .1cm

 We denote by $\Cal D^{\leq m}\Cal Man$
the category obtained from the category of dg-manifolds
by formally inverting all the $m$-quasiisomorphisms. 
 Isomorphism classes of objects in $\Cal D^{\leq m}\Cal Man$
 can be thus seen as analogous to  
$m$-homotopy types in topology.

\vskip .1cm

Returning to our particular situation, note that for any $p<q$
we have the morphisms of dg-manifolds
$$\alpha_{p,q}: RJ(k, A_{[p,q]})\to RJ(k, A_{[p+1, q]}),$$
$$\beta_{p,q}: RJ(k, A_{[p,q]})\to RJ(k, A_{[p, q-1]}),$$
given by the natural projection (forgetting one of the
graded components). These morphisms commute with each other
in the obvious sense. Now, the results of (4.1) can
be reformulated as follows.

\proclaim{(4.2.2) Proposition} Let $m\geq 0$
be given. Then there exists $p_0\geq 0$ such that for each $p\geq p_0$
 there exists $q_0$
such that for each $q\geq q_0$ the morphisms $\alpha_{p,q}$ and 
$\beta_{p,q}$ are $m$-quasiisomorphisms.
\endproclaim

\noindent {\sl Proof:} This follows from the equality (3.6.7),
Theorem 4.1.1 and Proposition 4.2.1. 

\proclaim {(4.2.3) Definition} The $m$-truncated derived Hilbert
scheme $RHilb_h^{\leq m}(X)$ is the object of $\Cal D^{\leq m}\Cal Man$
represented by any of the dg-manifolds $RJ(k, A_{[p,q]})$ where
$p, q$ are in the range given by Proposition 4.2.2. 
\endproclaim

Theorem 4.1.1 implies then the following

\proclaim{(4.2.4) Theorem} We have $\pi_0 RHilb_h^{\leq m}(X) = 
Hilb_h(X)$
and for each subscheme $Z\i X$ with sheaf of ideals $\Cal I\i \Cal O_X$ representing
a $\Bbb K$-point $[Z]\in Hilb_h(X)$,
$$H^iT^\bullet_{[Z]} RHilb_h^{\leq m}(X) = \op{Ext}^i_{\Cal O_Z}
(\Bbb L^\bullet_{Z/X}, \Cal O_Z), \quad 0\leq i\leq m.$$
\endproclaim

\vskip .2cm

\noindent {\bf (4.3) The case of locally complete intersections.}
In general, if $\Bbb L^\bullet_{Z/X}$ has infinitely many
cohomology sheaves, none of the above truncations
 is sufficient to capture all the
Ext's at once. We now concentrate on a special case when
such uniform truncation is possible. Recall that $X\i \Bbb P^n$
is assumed smooth. Consider the open subscheme
$Hilb_h^{LCI}(X)\i Hilb_h(X)$ formed by subschemes $Z\i X$ with
the Hilbert polynomial $h$,  for which
the embedding $Z\hookrightarrow X$ is 
a locally complete intersection morphism \cite{I}. 
As well known, for such $Z$ the relative
cotangent complex is quasiisomorphic to one sheaf
in degree $(-1)$ which is locally free and called the    conormal
bundle:
$$\Bbb L^\bullet_{Z/X} = N^*_{Z/X}[1].\leqno (4.3.1)$$ 
(In particular, if $X, Z$ are smooth, then we have (4.3.1)
with $N_{Z/X}^*$ being the conormal bundle in the usual sense.)
Therefore, for  such $Z$, 
$$\op{Ext}^i_{\Cal O_Z}(\Bbb L^\bullet_{Z/X}, \Cal O_Z) = H^{i-1}
(Z, N_{Z/X}), \leqno (4.3.2)$$
and this vanishes whenever $i\geq \dim(Z)$.

Recall that for any dg-scheme $X= (X^0, \Cal O^\bullet_X)$ and
any open subset $U\i X^0$ we have the induced dg-scheme
$X_U = (U, \Cal O^\bullet_X|_U)$. 

Let $0\ll p\ll q$ and let $U_{p,q}\i RJ(k, A_{[p,q]})^0$ be
the open subset obtained by removing the closed subset
in $Hilb_h(X) = \pi_0(RJ(k, A_{[p,q]})\buildrel \op{closed}\over
\hookrightarrow  RJ(k, A_{[p,q]})^0$ which is the complement
of $Hilb_h^{LCI}(X)$.
 The vanishing
of sufficiently high 
  Ext's in (4.3.2) 
 gives us the following.

\proclaim{(4.3.3) Proposition} For $0\ll p\ll q$ the dg-manifolds
$ RJ(k, A_{[p,q]})_{U_{p,q}}$ belong to the same quasiisomorphism
class (the quasiisomorphisms being established by the appropriate
restrictions of the morphisms $\alpha_{p,q}$ and $\beta_{p,q}$). 
  \endproclaim

We denote $RHilb_h^{LCI}(X)$ and call the {\it derived LCI-Hilbert
scheme} of $X$ the dg-manifold represented (up to
quasiisomorphism) by any of the $ RJ(k, A_{[p,q]})_{U_{p,q}}$ above.
If we need to emphasize the dependence of this manifold of $p$ and
$q$, we will use the notation $RHilb_h^{LCI}(X, [p,q])$. 
The following properties are now obvious from the construction.

\proclaim{(4.3.4) Theorem}
We have $\pi_0 RHilb_h^{LCI}(X) = Hilb_h^{LCI}(X)$ and for
any locally complete intersection $Z\i X$ representing
a $\Bbb K$-point $[Z]\in Hilb_h^{LCI}(X)$
$$H^i T^\bullet_{[Z]} RHilb_h^{LCI}(X) = H^{i} (Z, N_{Z/X}).$$
\endproclaim

\vskip .1cm

\noindent {\bf (4.3.5) Example: the derived space of maps.}
Let $C, Y$ be  projective schemes with $Y$ being smooth.
 The scheme $\op{Map}(C, Y)$
of maps $C\to Y$ can be seen as an open subscheme of the full
Hilbert scheme $Hilb(C\times Y)$ consisting of subschemes $C'\i C\times Y$
whose projection to $C$ is an isomorphism. This scheme is typically
an infinite disjoint union of schemes of finite type.
If we fix a projective embedding of $C\times Y$, then for any $h\in \Bbb Q[t]$
the intersection
 $$\op{Map}_h(C, Y) = \op{Map}(C, Y)\cap Hilb_h(C\times Y)$$
is a scheme of finite type which is open and closed
in $\op{Map}(C, Y)$. We can think of $h$ as being the ``generalized degree"
of a map $f: C\to Y$. Note that $\op{Map}_h(C, Y) \i Hilb_h^{LCI}
(C\times Y)$. 

Thus, deleting from $RHilb_h^{LCI}(X, [p,q])^0$ the closed
subset formed by points in $Hilb_h^{LCI}(X)$ which do not
represent graphs of maps, we
  get, by restriction, a dg-manifold $\op{RMap}_h(C, Y)$
with the following properties:
$$\pi_0 \op{RMap}_h(C, Y) = \op{Map}_h(C, Y); \leqno (4.3.5.1)$$
$$H^i T^\bullet_{[f]}\op{RMap}_h(C, Y) = H^i(C, f^*TY), \quad i\geq 0. \leqno
(4.3.5.2)$$
Here $f: C\to Y$ is a morphism representing a $\Bbb K$-point
$[f]$ of $\op{Map}_h(C, Y)$.  Again, this dg-manifold itself
depends, strictly speaking, on the choice of $0\ll p\ll q$ and
we will use the notation $\op{RMap}_h^{[p,q]}(C,Y)$
if we will need to emphasize this dependence. 

\vskip .3cm

\noindent {\bf (4.4) The relative version.} Let $S$ be a scheme
over $\Bbb K$ and $\pi: X\to S$ be a flat projective
morphism. We fix an $S$-embedding $X\hookrightarrow S\times \Bbb P^n$
and denote by $\Cal O(1)$ the corresponding  relative very ample sheaf.
We have then the sheaf of graded commutative algebras
$\Cal A = \bigoplus_i \Cal A_i$ with $\Cal A_i = \pi_* \Cal O(i)$
being coherent and locally free for $i\gg 0$. 
Thus for $0\ll p\ll q$  the dg-scheme $\bar RJ(k, \Cal A_{[p,q]}/S)$ 
from (3.7) is
defined and is smooth over $S$. We also have the relative versions of
the morphisms $\alpha_{p,q}$ and $\beta_{p,q}$ from (4.2)
which we denote by the same letters. 

\proclaim {(4.4.1) Proposition} Suppose $S$ is of finite type and let $m\geq 0$
be given. Then there exists $p_0\geq 0$ such that for each $p\geq p_0$
 there exists $q_0$
such that for each $q\geq q_0$ the relative
morphisms $\alpha_{p,q}$ and 
$\beta_{p,q}$ are $m$-quasiisomorphisms.
\endproclaim

\noindent {\sl Proof:}  As in the proof
of Proposition 4.2.2, we use the equality (3.6.7),
Theorem 4.1.1 and Proposition 4.2.1. In addition, we use
Proposition 3.7.1 to get the information about
the relative 
tangent dg-spaces to $\bar RJ(k, \Cal A_{[p,q]}/S)$
at arbitrary $\Bbb K$-points. 

\vskip .2cm

We define now the relative $m$-truncated derived
Hilbert scheme $RHilb_h^{\leq m}(X/S)$ to be
the object of $\Cal D^{\leq m}\Cal Man$
represented by any of the dg-schemes $\bar RJ(k, \Cal A_{[p,q]}/S)$ where
$p, q$ are in the range given by Proposition 4.4.1.
Further, the definition of the relative derived LCI-Hilbert
scheme $RHilb_h^{LCI}(X/S)$ is completely analogous to
that given in (4.3) and is left to the reader. By
construction and by Theorem 1.4.1 we have
$$\pi_0 RHilb_h^{\leq m}(X/S) = Hilb_h(X/S), \quad\pi_0 RHilb_h^{LCI}(X/S) =
Hilb_h^{LCI}(X/S)$$
as well as a compatibility statement analogous to Proposition 3.7.1(b).

We will use the notation $RHilb^{LCI}(X/S, [p,q])$ if we need to
emphasize the dependence on the choice of $p,q$. 

\vskip .1cm

Let $S$ be a scheme of finite type,
and $C, Y$ be flat projective schemes over $S$ with $Y$ being smooth.
Then we have the relative derived space of maps
$\op{RMap}_h(C, Y|S)$ 
(or $\op{RMap}_h^{[p,q]} (C, Y|S)$, if we need
to emphasize the dependence on the choice of $0\ll p\ll q$)
defined, as in Example 4.3.6,
via the relative derived LCI-Hilbert scheme $RHilb_h^{LCI}(C\times_S Y,
[p,q])$.
The following fact which we record for future reference,
follows directly from the constructions and Proposition 3.7.1.

\proclaim{(4.4.2) Proposition} Let  $\phi: S'\to S$
be a
 morphism of  schemes of finite type  and $C_{S'}, Y_{S'}$ be the 
pullbacks of $C, Y$ to schemes over $S'$.  If $0\ll p\ll q$
are admissible for $S, C, Y$, then they
are admissible for $S', C_{S'}, Y_{S'}$ as well and
$\op{RMap}_h^{[p,q]}(C_{S'}, Y_{S'}|S')$ is isomorphic to the
 fiber product of $S'$ and $\op{RMap}_h^{[p,q]}(C, Y|S)$ over $S$. 
\endproclaim

\vfill\eject

\heading { 5. The relative version over a stack and
derived moduli stacks of stable maps}
\endheading

In this section we  further relativize the above
constructions to the case when the base scheme $S$ in (3.7) and (4.4)
is replaced by an algebraic stack of finite type. 
As an application, we construct the derived
stacks of stable maps, thus realizing the suggestion of
Kontsevich \cite{Kon}. 

Our basic reference for stacks is \cite{LM}
from which we adopt all the terminology and conventions. 
As before, we  work over a field $\Bbb K$ of characteristic 0.

\vskip .2cm

\noindent {\bf (5.1) Dg-stacks.}  We first develop some minimal formalism
necessary to speak about derived moduli stacks. While certainly
too restrictive for all applications (e.g., for the construction
of the derived moduli stack of vector bundles), this formalism
will be sufficient for our problem.

\proclaim{(5.1.1) Definition} A dg-stack is a pair
$\Cal X = (\Cal X^0, \Cal O^\bullet_{\Cal X})$ where $\Cal X^0$
is an algebraic stack and $\Cal O^\bullet_{\Cal X}$ is a $\Bbb Z_-$-graded
quasicoherent $\Cal O_{\Cal X^0}$-dg-Algebra such that
$\Cal O_{\Cal X}^0 = \Cal O_{\Cal X^0}$. A graded stack is a dg-stack
with trivial differential.  1- and 2- morphisms of dg-stacks
are defined as 1- and 2-morphisms of (dg-)ringed stacks
(\cite{LM}, Def. 12.7.2). 

\endproclaim 

We denote by $St$ and $dgSt$ the 2-categories of algebraic
stacks and dg-stacks respectively. 

The $\Cal O_{\Cal X^0}$-Algebra $\underline{H}^0 (\Cal O^\bullet_{\Cal X})$
is a quotient of $\Cal O_{\Cal X^0}$ and so defines a closed
substack
$$\tau_{\leq 0}\Cal X = \op{Spec}(\underline{H}^0(\Cal O^\bullet_{\Cal X}))\i
\Cal X^0.$$
See \cite{LM}, (14.2.7) for the meaning of Spec in this situation. 
It is clear that for any algebraic stack $\Cal Y$ we have an equivalence
of categories
$$\op{Hom}_{dgSt}(\Cal Y, \Cal X) \simeq \op{Hom}_{St}(\Cal Y, \tau_{\leq 0}
\Cal X).$$
Further, the sheaf of graded commutative
algebras $\underline{H}^\bullet(\Cal O_{\Cal X})$ is naturally
a sheaf on $\tau_{\leq 0} \Cal X$. We denote by
$\Cal X_h$ the graded stack
$(\tau_{\leq 0} \Cal X, \underline{H}^\bullet(\Cal O_{\Cal X}))$.
A 1-morphism of dg-stacks $F: \Cal X\to\Cal Y$ is called a quasi-equivalence
if $F_h: \Cal X_h\to\Cal Y_h$ 
is an equivalence of ringed stacks. 
The concept of an $m$-quasi-equivalence is defined similarly,
by considering the truncations of $\underline{H}^\bullet
(\Cal O_{\Cal X}^\bullet)$ and $\underline{H}^\bullet
(\Cal O_{\Cal Y}^\bullet)$ in degrees $\geq -m$. 

We will also denote by $\Cal X_\sharp$ the graded stack
$(\Cal X^0, \Cal O_{\Cal X, \sharp}^\bullet)$,
where $\Cal O_{\Cal X, \sharp}^\bullet$ is $\Cal O_{\Cal X}^\bullet$
with forgotten differential.

\proclaim{(5.1.2) Definition}
A dg-stack $\Cal X$ is said to be of (locally) finite type,
if $\Cal X^0$ is an algebraic stack of (locally) finite type
and each $\Cal O_{\Cal X}^i$ is a coherent $\Cal O_{\Cal X^0}$-Module.
\endproclaim

\proclaim{(5.1.3) Definition}
Let  $F: \Cal X\to\Cal Y$ be a 1-morphism of dg-stacks of finite type.
We say that $F$ is smooth, if for any presentation 
(in the sense of \cite{LM}, Def. 4.14) of the morphism $F^0: \Cal X^0
\to\Cal Y^0$
$$\matrix
X'&\buildrel \beta\over\lra &\Cal X^{0'}&\buildrel\alpha\over\lra& \Cal X^0&\\
&f\searrow &\big\downarrow& \square&\big\downarrow&F^0\\
&&Y&\buildrel\gamma\over\lra &\Cal Y^0&\\
\endmatrix 
\leqno (5.1.4)$$
with $X', Y$ schemes of finite type, the morphism of dg-schemes of finite
type
 $$(X', \beta^*\alpha^*\Cal O^\bullet_{\Cal X})\to (Y, \gamma^*\Cal O^\bullet
_{\Cal Y})$$
 is smooth in the sense of \cite{CK}, Def. 2.7.1.

\endproclaim

As usual, we say that $\Cal X$ is smooth if $\Cal X\to \{pt\}$ is a smooth
morphism. Thus, in our terminology, every smooth dg-stack
is of finite type. 

\proclaim{(5.1.5) Definition}
Let $\Cal X$ be a smooth dg-stack, $\Bbb F$ be a field extension of
$\Bbb K$  and $x$ be an $\Bbb F$-point of $\Cal X$
(i.e., in fact an $\Bbb F$-point of $\tau_{\leq 0}\Cal X$). The tangent
dg-space $T^\bullet_x \Cal X$ is a complex of $\Bbb F$-vector spaces situated
in degrees $[-1, \infty)$ which is defined (up to a unique isomorphism
in the derived category) as follows. Let $P: X^0\to\Cal X^0$
be a presentation of $\Cal X^0$ with $X^0$ being a smooth 
algebraic variety and let $X$ be the dg-manifold
$(X^0, P^*\Cal O^\bullet_{\Cal X})$. Then we set
$$T^\bullet_x\Cal X = \biggl\{ (T_{X^0/\Cal X^0})_{x'} \to T^0_{x'}X \to
T^1_{x'}X\to ...\biggr\},$$
where $x'
\in X^0$ is any $\Bbb K$-point over $x$ and  $(T_{X^0/\Cal X^0})_{x'}$
 denotes the fiber at $x'$ of the relative tangent bundle of
$X^0$ over $\Cal X^0$, which is embedded
into $T^0_{x'}X = T_{x'}X^0$ in the standard way.
\endproclaim

As in the case of dg-manifolds \cite{CK}, we denote
$$\pi_i(\Cal X, x) = H^i T^\bullet_x\Cal X, \quad i=1,0,-1, -2, ...$$
As in that case, we see easily that:

\proclaim{(5.1.6) Proposition}
We have well-defined ``Whitehead products"
$$[-,-]: \pi_i(\Cal X, x)\otimes\pi_j(\Cal X, x) \to \pi_{i+j-1}(\Cal X, x),$$
making $\pi_{\bullet -1}(\Cal X, x)$ into a graded Lie algebra. 
In particular, $\pi_1(\Cal X, x)$ is a Lie algebra in the ordinary sense,
and it is identified with the Lie algebra of the algebraic group
$\op{Aut}(x)$. 
\endproclaim

\proclaim{(5.1.7) Proposition} Let $F: \Cal X\to\Cal Y$ be a 1-morphism
of smooth dg-stacks. Then the following are equivalent:

(i) $F$ is a quasi-equivalence (resp. an $m$-quasi-equivalence, $m\geq 0$).

(ii) The 1-morphism of algebraic stacks $\tau_{\leq 0}F: \tau_{\leq 0}\Cal X
\to\tau_{\leq 0}\Cal Y$ is an equivalence and for any field extension
$\Bbb F\supset \Bbb K$ and any
$\Bbb F$-point $x$ of $\Cal X$ the  natural morphism of complexes
$d_xF: T^\bullet_x\Cal X\to T^\bullet_{F(x)}\Cal Y$ is
a quasi-isomorphism (resp. an $m$-quasi-isomorphism). 

\endproclaim

\noindent {\sl Proof:}
This follows from the analogous
result for dg-manifolds (\cite{CK}, Prop. 2.5.9) and the following
obvious remark applied to the morphism of sheaves
$(\tau_{\leq 0} F)^*\underline{H}^\bullet(\Cal O_{\Cal Y}^\bullet)\to
\underline{H}^\bullet(\Cal O_{\Cal X})$. 

\proclaim{(5.1.8) Lemma} Let $\Cal Z$ be an algebraic stack
(of finite type)
and let $f: \Cal F\to\Cal G$ be a morphism of quasi-coherent
sheaves on $\Cal Z$. Then $f$ is an isomorphism if and only if for any
presentation $P: Z\to \Cal Z$ (with $Z$ of finite type)
the induced morphism $P^*\Cal F
\to P^*\Cal G$ is an isomorphism. 
 
\endproclaim

\vskip .2cm

\noindent {\bf (5.2) Dg-schematic local construction.} We recall
the formalism of \cite{LM}, \S 14, slightly modified for our
purposes. 

\proclaim{(5.2.1) Definition} Let $\Cal X$ be an algebraic stack.
A dg-schematic local construction on $\Cal X$ is a datum $\underline{\Cal Y}$
which associates to any morphism $U\to X$ where $U$ is a scheme,
a dg-scheme $\underline{\Cal Y}_U$ and, further, to any morphism $\phi: V\to U$
of schemes, associates an isomorphism of dg-schemes 
$V\times_U \underline{\Cal Y}_U$ so that these isomorphisms satisfy the
cocycle conditions. 

\endproclaim

\proclaim {(5.2.2) Proposition} (a) Any dg-schematic local
construction $\underline{\Cal Y}$ on $\Cal X$ gives rise to a dg-stack
$\Cal Y$ equipped with a morphism $\Cal Y\to\Cal X$ so that
for any $U\to \Cal X$ as above the dg-stack $U\times_{\Cal X}\Cal Y$
is equivalent to the dg-scheme $\underline{\Cal Y}_U$.

(b) Further, the algebraic stack $\tau_{\leq 0}\Cal Y$ corresponds
to the local construction $U\mapsto \pi_0 \underline{\Cal Y}_U$.
\endproclaim

\noindent {\sl Proof:} (a) When the dg-structures on the $\underline{\Cal Y}_U$
are trivial, this is Proposition 14.1.7 of \cite{LM} (with algebraic
spaces replaced by schemes). One proves, in a similar way,
a version of the cited statement in which the $\underline{\Cal Y}_U$
are not just schemes but ringed schemes (i.e., schemes with
a quasicoherent sheaves of rings) and, further, dg-ringed schemes.
This implies (a). Part (b) is equally clear. 

\vskip .2cm

If $\Cal X$ is a stack of finite type, then in the above
discussion it is clearly sufficient to consider the schemes $U$
of finite type, which we will do. 

\vskip .3cm

\noindent {\bf (5.3) Derived spaces of maps over a stack.}
Let $\Cal S$ be an algebraic stack of finite type and
$\Cal A$ be a commutative $\Cal O_{\Cal S}$-Algebra which, as an
 $\Cal O_{\Cal S}$-Module, is locally free of finite rank. Fix $k\geq 0$.
Application of the concept of a schematic local construction
gives the relative Grassmannian
$G(k, \Cal A/\Cal S)$ and the relative ideal space
$J(k, \Cal A/\Cal S)$ which are algebraic stacks equipped with
schematic projective morphisms to $\Cal S$ and 
$J(k, \Cal A/\Cal S)$ being a closed substack in 
$G(k, \Cal A/\Cal S)$. 

Notice further that associating to any morphism $\phi: U\to \Cal S$ with
$U$ a scheme, the dg-scheme $\bar RJ(k, (\phi^*\Cal A)/U)$, see (3.7),
is a dg-schematic local construction in virtue of Proposition 3.7.1(b).
Thus we get a dg-stack $\bar RJ(k, \Cal A/\Cal S)$ with
$\tau_{\leq 0} \bar RJ(k, \Cal A/\Cal S) = J(k, \Cal A/\Cal S)$.

Similarly, if $\Cal A = \bigoplus_{i=a}^b \Cal A_i$ is a graded
$\Cal O_{\Cal S}$-Algebra with each $\Cal A_i$ being
locally free of finite rank as an $\Cal O_{\Cal S}$-Module
and $k= (k_i)_{i=a}^b$, we get the dg-stack of graded ideals
which we still denote $\bar RJ(k, \Cal A/\Cal S)$.

 \vskip .2cm

Next, assume that $\Cal S$ is of finite type and
let $\pi: \Cal X\to\Cal S$ be a flat schematic projective morphism
of algebraic stacks and $\Cal O(1)$ be the fixed relatively very ample
sheaf on $\Cal X$. Then setting $\Cal A_i = \pi_*\Cal O(i)$,
we find that for $i\gg 0$ each $\Cal A_i$ is locally free of finite rank
as an $\Cal O_{\Cal S}$-Module and so $\Cal A_{[p,q]} = \bigoplus_{i=p}^q
\Cal A_i$ is an $\Cal O_{\Cal S}$-Algebra satisfying the above conditions,
so  the dg-stack $\bar RJ(k, \Cal A_{[p,q]}/\Cal S)$
is defined. Because
$\Cal S$ is of finite type, we see, by using Lemma 5.1.8,
that for $0\ll p\ll q$ the dg-stacks 
$\bar RJ(k, \Cal A_{[p,q]}/\Cal S)$ have the same $m$-quasi-equivalence
class and thus can be, with some abuse of notation,
denoted by $RHilb_h^{\leq m}(\Cal X/\Cal S)$. By construction,
$$\tau_{\leq 0} RHilb_h^{\leq m}(\Cal X/\Cal S) = 
 Hilb_h(\Cal X/\Cal S) \leqno (5.3.1)$$
is the relative Hilbert scheme of $\Cal X$ over $\Cal S$. 
Further, associating to any $\phi: U\to \Cal S$ the dg-scheme
$RHilb_h^{LCI}((\phi^*\Cal X/U, [p,q])$, see (4.2.4),
forms again a dg-schematic construction so we get
a dg-stack $RHilb_h^{LCI}(\Cal X/\Cal S, [p,q])$
whose quasi-equivalence type for $0\ll p\ll q$ is
independent on the choice of such $p,q$. 

\vskip .2cm

Finally, let $\Cal C\to \Cal S$, $\Cal Y\to\Cal S$ be flat
schematic projective morphisms, with $Y$ being smooth.
By the same token as before, we get
the relative space of maps $RMap_h^{[p,q]}(\Cal C, \Cal Y|\Cal S)$
whose  quasi-equivalence type for $0\ll p\ll q$ is
independent on the choice of such $p,q$.

\vskip .3cm

\noindent {\bf (5.4) The derived stack of  stable maps.} 
 Let $\widetilde{M}_{g, n}$  be the stack
of prestable $n$-pointed curves of genus $g$, see \cite{BM}, \S 2.
Thus an $S$-point of this stack is a system
$(\pi: C\to S, x_1, ..., x_n)$ where $\pi$ is a
 flat proper
morphism $\pi: C\to S$ of relative dimension 1 such that every geometric
fiber $C_s, s\in S$, is a connected projective curve of arithmetic genus $g$,
with at most nodal singularities and 
$x_i: S\to C$ are sections
whose  values at every geometric point $s\in S$
are smooth and distinct points of the curve $C_s$.
We denote by $\Cal C\to \widetilde{M}_{g, n}$ the universal
curve. This is a schematic projective 1-morphism
of stacks. 

If $K$ is a field and $(C, x_1, ..., x_n)$ is an $n$-pointed
prestable curve of genus $g$ over $K$,
 we denote by $[C, x_1, ..., x_n]$ the corresponding $K$-point
of $\widetilde{M}_{g, n}$. As for any $K$-point
of any algebraic stack, the tangent space $T^\bullet_{[C, x_1, ...,
x_n]}\widetilde {M}_{g, n}$ is a complex of $K$-vector spaces,
situated in degrees $[-1, 0]$ and defined up to quasiisomorphism.
The first order deformation theory
gives an identification
$$T^\bullet_{[C, x_1, ..., x_n]}\widetilde{M}_{g, n} = R\Gamma(C, RTC(-x_1-...-x_n))[1],
\leqno (5.4.1)$$
where $RTC$ is the tangent complex of $C$. 

\proclaim{(5.4.2) Definition}
Let $(Y, \Cal O(1))$ be a smooth projective variety over $\Bbb K$,
let $d\in \Bbb Z$ 
and $S$ be a $\Bbb K$-scheme. An $n$-pointed prestable
map of genus $g$ and degree $d$ over $S$ is a system
$(C, x_1, ..., x_n, f)$ where $\pi: C\to S$ and $x_i: S\to C$ are
as in (5.1.1) while $f: C\to S\times Y$ is an $S$-map
such that for any geometric point $s\in S$ the degree of the sheaf
$f^*\Cal O(1)$ on $C_s$ is $d$. 
\endproclaim

 We denote by $\widetilde{M}_{g, n}(Y, d)$ the stack of
prestable maps, as in (5.1.2). This is an algebraic stack which
is non-separated and possibly non-smooth. By definition,
we have the following, cf. \cite{B}, p. 604.

\proclaim{(5.4.3) Proposition} The stack $\widetilde{M}_{g, n}(Y, d)$
is naturally identified with the relative space of maps
 $\op{Map}_h(\Cal C, \widetilde{M}_{g, n}\times Y|\widetilde{M}_{g, n})$
where $h\in \Bbb Q[t]$ is a polynomial uniquely determined by $g, d$
and the sheaf $\Cal O(1)$ on $Y$.
\endproclaim

As before, if $K$ is a field
and $(C, x_1, ..., x_n, f)$ is an $n$-pointed prestable map
over $K$, we denote by $[C, x_1, ..., x_n, f]$ the corresponding
$K$-point of $\widetilde{M}_{g,n}(Y, d)$. The first
order deformation theory now identifies
$$T^\bullet_{[C, x_1, ..., x_n, f]} \widetilde{M}_{g,n}(Y, d) = 
\leqno (5.4.4)$$
$$\biggl( \tau_{\leq 1} R\Gamma \biggl(C, \op{Cone}
\bigl\{ RTC(-x_1-...-x_n)\buildrel df\over\lra f^*TY\bigr\}\biggr)\biggr) [1],$$
where $\tau_{\leq 1}$ is the cohomological truncation of a complex
in degrees $\leq 1$. 
The $(-1)$st cohomology of (5.4.4) is (as is the case with any 
algebraic stack) 
the Lie algebra of infinitesimal automorphisms of $(C, x_1, ..., x_n, f)$.

\proclaim{(5.4.5) Definition} A prestable map $(C, x_1, ..., x_n, f)$
over a field $\Bbb F$ is called stable, if
 $H^{-1} T^\bullet_{[C, x_1, ..., x_n, f]} \widetilde{M}_{g,n}(Y, d) = 0$.
A prestable map over a $\Bbb K$-scheme $S$ is called stable, if for any 
geometric point $s\in S$ the induced map over the field $\Bbb K(s)$
is stable. 
\endproclaim

We denote by $\overline{M}_{g,n}(Y, d)\i \widetilde{M}_{g,n}(Y,d)$
the open substack formed by stable maps. It is known \cite{Kon}\cite{BM} that
$\overline{M}_{g,n}(Y, d)$ is a proper Deligne-Mumford
stack, in particular, it is separated and  of finite type.
 However, it is not smooth in general,
essentially because its tangent spaces, given by the same
formula (5.4.4), are obtained by truncating some naturally arising
complex.  We now proceed to construct a smooth derived version of 
$\overline{M}_{g,n}(Y, d)$.

\vskip .1cm

Forgetting the map $f$ defines  morphisms of stacks
$$\overline{\rho}: \overline{M}_{g,n}(Y, d)\to \widetilde{M}_{g,n},
\quad \widetilde{\rho}: \widetilde {M}_{g,n}(Y, d)\to \widetilde{M}_{g,n} .
\leqno (5.4.6)$$
Note that though $\widetilde{M}_{g,n}$ is  only locally of finite type, $\overline{M}_{g,n}(Y, d)$
is of finite type, so the image of $\overline{\rho}$ is contained
in an open substack  $\Cal S$ of  $\widetilde{M}_{g,n}$ of finite type.

\proclaim {(5.4.7) Definition} Let $\Cal S\i \widetilde{M}_{g,n}$ be
an open substack of finite type. The derived stack of
prestable maps of type $\Cal S$ is defined to be
$$R\widetilde{M}_{g,n}^{\Cal S}(Y,d) = \op{RMap}_h(\Cal C, Y\times\Cal S |
\Cal S),$$
where $h$ is as in (5.4.3).
 
\endproclaim

\proclaim{(5.4.8) Theorem}
(a) $R\widetilde{M}^{\Cal S}_{g,n}(Y,d)$ is a smooth dg-stack
and  $\tau_{\leq 0} R\overline M_{g,n}^{\Cal S}(Y,d) =
 \widetilde{\rho}^{-1}(\Cal S)$.

(b) If  $(C, x_1, ..., x_n, f)$ is an $n$-pointed prestable map
over $\Bbb K$ such that $(C, x_1, ..., x_n)\in\Cal S$,
then
$$T^\bullet_{[C, x_1, ..., x_n, f]} R\widetilde{M}_{g,n}^{\Cal S}(Y, d) = $$
$$ = R\Gamma \biggl(C, \op{Cone}
\bigl\{ RTC(-x_1-...-x_n)\buildrel df\over\lra f^*TY\bigr\}\biggr) [1].$$
\endproclaim

\noindent {\sl Proof:} (a) Smoothness of 
$R\widetilde{M}_{g,n}^{\Cal S}(Y, d)$ follows from the smoothness
of $\Cal S$ and from the relative smoothness of RMap. The identification
of $\tau_{\leq 0}$ follows from Proposition 5.4.3 and from the
identification $\pi_0\op{RMap} = \op{Map}$. 

(b) Let us write $R\widetilde{M}^{\Cal S}$ for 
$ R\widetilde{M}_{g,n}^{\Cal S}(Y, d)$ and $x$ for $x_1, ..., x_n$.
Then we have the distinguished triangle of complexes (relative
vs. asolute tangent dg-spaces)
$$T^\bullet_{[C, x, f]} (R\widetilde{M}^{\Cal S}/\Cal S)
\to
T^\bullet_{[C, x, f]} R\widetilde{M}^{\Cal S} \to T^\bullet_{[C, x]}\Cal S.$$
Now, $T^\bullet_{[C, x]}\Cal S$ is given by (5.4.1), while by Proposition
4.4.2
$$T^\bullet_{[C, x, f]} (R\widetilde{M}^{\Cal S}/\Cal S) = 
T^\bullet_{[f]} \op{RMap}_h(C, Y) = R\Gamma(C, f^*TY)$$
and our assertion follows from comparing the above triangle with the
triangle
$$R\Gamma(C, f^*TY)\to R\Gamma(C, \op{Cone}) \to R\Gamma(C, RTC(-x_1-...-x_n))
[1],$$
where Cone is the cone of $df$ in the formulation of (b). Theorem is
proved.

\vskip .2cm

For a dg-stack $\Cal X$ and an open substack $\Cal U\i\Cal X^0$
we have a dg-stack $\Cal X|_{\Cal U} = (\Cal U, \Cal O^\bullet_{\Cal X}|_{\Cal
U})$. Further, $\tau_{\leq 0}\Cal X$ being a closed substack of $\Cal X^0$,
any open subset $\Cal V\i \tau_{\leq 0}\Cal X$ gives rise to
an open substack $\Cal V'\i\Cal X^0$
which is the complement in $\Cal X^0$ to
the complement of $\Cal V$ in $\tau_{\leq 0} \Cal X$ considered as 
a closed substack in $\Cal X$. Let us write $\Cal X_{\Cal V}$
for $\Cal X|_{\Cal V'}$.

\proclaim{(5.4.9) Definition}
The derived moduli stack of stable maps is defined to be
$$R\overline{M}_{g,n}(Y,d) = R\overline{M}_{g,n}^{\Cal S}(Y,d)_{
\overline{M}_{g,n}(Y,d)},$$
where $\Cal S\i \overline{M}_{g,n}$ is any open substack
containing the image of $\overline{\rho}$ in (5.4.6).
\endproclaim

Definition 5.4.9 gives, in particular, a construction of the sheaf
of graded commutative algebras $\underline{H}^\bullet(\Cal O^\bullet_{
R\overline{M}_{g,n}(Y,d)})$ on the Deligne-Mumford stack
$\overline{M}_{g,n}(Y,d)$. It is this structure whose existence
was originally conjectured by Kontsevich (\cite{Kon}, n. 1.4.2)
and used to give a formula for the virtual fundamental class
({\it ibid.} p. 344). In a subsequent paper we plan
to study this fundamental class in more detail
and, in particular, show that it conincides with the class
constructed in \cite{B} \cite{BF} \cite{LT}.

\vfill\eject

\enddocument